\documentclass[11pt,bezier]{article}  
\usepackage[textsize=tiny]{todonotes}
\usepackage{color}              
\usepackage{epsfig}
\usepackage{graphicx}
\usepackage{amsmath}
\usepackage{amsfonts}
\usepackage{amsthm}
\usepackage{amssymb}
\usepackage{epsfig}
\usepackage{hyperref}
\usepackage{bbm}
\usepackage{longtable}

\definecolor{darkgreen}{rgb}{0,.3,0}
\definecolor{darkagenta}{rgb}{.5,0,.5}
\definecolor{darkred}{rgb}{1,0,0}
\definecolor{darkblue}{rgb}{0,0,.4}

\newcommand{\ch}[1]{#1}
\newcommand{\chr}[1]{#1}
\newcommand{\chn}[1]{#1}
\newcommand{\chrr}[1]{#1}

\renewcommand{\P}{\mathbb{P}}

\renewcommand{\P}{\mathbb{P}}
\newcommand{\E}{\mathbb{E}}

\def\N{{\mathbb N}}
\def\R{{\mathbb R}}

\def\P{{\mathbb P}}
\def\E{{\mathbb E}}

\newcommand{\e}{{\mathrm e}}

\newcommand{\eqn}[1]{\begin{equation} #1 \end{equation}}
\newcommand{\eqan}[1]{\begin{align} #1 \end{align}}
\newcommand{\lbeq}[1]{\label{#1}}
\newcommand{\refeq}[1]{(\ref{#1})}
\newcommand{\sss}{\scriptscriptstyle}
\newcommand{\Var}{{\rm Var}}
\newcommand{\Cov}{{\rm Cov}}
\newcommand{\expec}{{\mathbb E}}
\newcommand{\prob}{{\mathbb P}}

\newcommand{\op}{o_{\sss \prob}}

\newcommand{\Op}{O_{\sss \prob}}

\newcommand {\vep}{\varepsilon}
\newcommand {\convd}{\stackrel{d}{\longrightarrow}}
\newcommand {\convp}{\stackrel{\sss {\mathbb P}}{\longrightarrow}}

\newcommand {\convas}{\stackrel{a.s.}{\longrightarrow}}
\newcommand{\nn}{\nonumber}

\newcommand{\ra}{\rightarrow}

\newcommand\1{\mathbbm{1}}
\newcommand{\indic}[1]{\1_{\{#1\}}}

\numberwithin{equation}{section}

\newtheorem{Lemma}{Lemma}[section]

\newtheorem{Theorem}[Lemma]{Theorem}

\newtheorem{Condition}[Lemma]{Condition}

\begin{document}

\title{Degree-degree dependencies in random graphs\\
with heavy-tailed degrees
}

\author{Remco van der Hofstad\footnote{Eindhoven University of Tedhnology and Eurandom} and Nelly Litvak\footnote{University of Twente}}

\maketitle

\begin{abstract}
Mixing patterns in large self-organizing networks, such as the Internet, the World Wide
Web, social and biological networks are often characterized by degree-degree dependencies between neighbouring nodes. In {\it assortative} networks\chrr{,} the degree-degree dependencies are positive (nodes with similar degrees tend \chrr{to} connect to each other), while in {\it disassortative} networks\chrr{,} these dependencies are negative. One of the problems with the commonly used \chn{Pearson correlation coefficient, also known as the {\it assortativity coefficient}} is that its magnitude decreases with the network size \chr{in disassortative networks}. This makes it impossible to compare mixing patterns, for example, in two web crawls of different sizes. As an alternative, we have recently suggested to use rank correlation measures, such as Spearman's rho. Numerical experiments have confirmed that Spearman's rho produces consistent values in graphs of different sizes but similar structure, and it is able to reveal strong (positive or negative) {dependencies} in large graphs.

In this paper we analytically investigate degree-degree dependencies for scale-free graph sequences. 
\chn{In order to demonstrate the ill behaviour of the Pearson's correlation coefficient,} we first study a simple model of two heavy-tailed highly correlated random variables $X$ and $Y$, and show that the sample correlation coefficient converges in distribution either to a proper random variable on $[-1,1]$, or to zero, \chrr{and the limit is non-negative a.s.\ if $X,Y\ge 0$.}
We next adapt these results to the \chn{degree-degree dependencies} in networks
as described by the \chn{Pearson} correlation coefficient, 
and show that it is non-negative in the large graph limit 
when the \chr{asymptotic} degree distribution has an infinite third moment. Furthermore, we provide examples where the Pearson's correlation coefficient converges to zero in a network with strong negative degree-degree dependencies, and another example where this coefficient converges in distribution to a random variable. We suggest the alternative degree-degree dependency measure, based on Spearman's rho, and prove that this statistical estimator converges to an appropriate limit under \chr{quite} general conditions. These conditions are proved to hold in common network models, such as \chr{the} configuration model and \chr{the} preferential attachment model. We conclude that rank correlations provide a suitable and informative method for uncovering network mixing patterns. 
\end{abstract}


\section{Introduction}
\label{sec:intro}

In this paper we present an analytical study of degree-degree correlations in graphs with power law degree distribution. In simple words, a random variable $X$ has a power-law distribution
with tail exponent $\gamma>0$ if its tail probability $\P(X>x)$ is
roughly proportional to $x^{-\gamma}$, for large enough $x$. Large self-organizing networks, such as the Internet, the World Wide Web, social and biological networks, usually exhibit \ch{high variation in the values of the degrees. Such networks are called \emph{scale free} \chr{indicating that there is no typical scale for the degrees}, and the \chr{high} \chn{degree} vertices are called \emph{hubs}.  This phenomenon is often modelled by using power-law degree distributions.}

Power-law distributions are {\it heavy tailed} since the tail
probability decreases much more slowly than \chr{a} negative exponential, and
thus one observes extremely large values of $X$ much more frequently
than in the case of light tails.  Statistical analysis of \ch{scale-free complex networks} 
has received massive attention in recent literature, see
e.g.~\cite{Mitzenmacher2004,Newman03survey}
for excellent surveys. Nevertheless, there still are many fundamental open problems. 
One of them is how to measure \emph{dependencies} between network
parameters.

An important \chn{characteristic} of networks is the dependency between
the degrees of direct neighbours. \chn{A network is usually called {\it assortative} when nodes with similar degrees are often connected, thus, the degree-degree dependencies are positive, while in a {\it disassortative} network these dependencies are negative.} \ch{The {\it degree-degree} dependencies define many  of the network's properties. For instance, the negative degree-degree correlations in the Internet graph have a great influence on the robustness to failures~\cite{Doyle2005ryf}, efficiency of Internet protocols~\cite{Li2005willinger}, {as well as} distances and betweenness~\cite{Mahadevan2006Internet_topology}. The correlation between in- and out-degree of tasks plays and important role in the  dynamics of production and development systems~\cite{Braha2007communication}. Mixing patterns affect epidemic spread~\cite{Eguiluz2002epidemic,Eubank2004nature} and Web ranking~\cite{Fortunato2007mean_field}.} 

Often, degree-degree dependence is
characterized by the \emph{assortativity coefficient} of the network,
introduced by Newman in \cite{Newman2002assortative}. The assortativity coefficient
is \chn{in fact \chrr{the} Pearson} correlation coefficient between the vector of degrees on each 
side of an edge, as a function of all edges. 
See \cite[Table I]{Newman2002assortative} for a list of assortativity coefficients
for various real-world networks. The empirical data suggest that
social networks tend to be assortative (the assortativity coefficient is
positive), while Internet, World Wide Web, and biological networks tend to be
disassortative. In \cite[Table I]{Newman2002assortative}, it is striking that, 
typically, larger disassortative networks have an 
assortativity coefficient that is closer to 0 and therefore
appear to have approximate {\it uncorrelated} degrees
across edges.  Similar conclusions can be drawn from \cite{Newman2003mixing},
see in particular \cite[Table II]{Newman2003mixing}. This phenomenon arises because Pearson's correlation coefficient in scale-free networks with realistic parameters decreases with the network size, as was pointed out in several recent papers~\cite{Dorogovtsev2010correlationsPA,Raschke2010correlations,HofLit13}. In this paper, we \chn{prove} that \chn{Pearson's} \chrr{correlation} coefficient in scale-free networks shows several types of pathological behavior, in particular, its infinite volume limit, \chr{when it exists,} is non-negative, independently of the mixing pattern, and in fact this limit can even be \emph{random}. 

In \cite{HofLit13} we propose an alternative measure for the degree-degree dependencies, based on the \emph{ranks} of degrees. This rank correlation 
approach is in fact classical in multivariate analysis, 
falling under the category of `concordance measures' - dependency measures 
based on \emph{order} rather than exact values of two stochastic variables. 
The huge advantage of such dependency measures is that they
work well \emph{independently} of the number of finite moments
of the degrees,  while \chn{Pearson's} coefficient suffers from a strong dependence on the extreme values of
the degrees. \chr{Recent applications of rank correlation measures, such as Spearman's rho~\cite{Spearman1904} and the closely related Kendall's tau~\cite{Kendall1938}, include the} concordance between two rankings for a set of documents in web search. In this application field many other measures for rank distances have been proposed, see e.g.~\cite{Kumar2010rank_correlations} and the references therein.

We show mathematically that statistical estimators for degree-degree dependencies based on rank correlations are \emph{consistent}. That is, for graphs of different sizes but similar structure (e.g.\ preferential attachment graphs of increasing size), these estimators converge to their `true' \chr{or limiting} value that describes the degree-degree dependence in an infinitely large graph (in particular, the variance of the estimator decreases as the size of the graph grows). We also show that \chn{Pearson's} correlation coefficient 
does not have this basic property when degree distributions are heavy-tailed. 
In particular, as explained in more detail in \cite{HofLit13}, this implies that the assortativity 
coefficient \chn{as suggested in \cite{Newman2002assortative}} does not allow one to compare the \chn{degree-degree dependencies in} graphs of different sizes, such as 
they arise when studying a network at different time stamps, or comparing two different networks, e.g.\, web crawls of different domains or Wikipedia graphs from different languages. On the other hand, such a comparison
\emph{is} possible using Spearman's rho. This paper forms the mathematical justification of our paper \cite{HofLit13}, where 
similar results were predicted on a less formal level and confirmed by numerical experiments.


\chn{The paper is organized as follows. In Section~\ref{sec-corr} we start with the analysis of the sample \chn{Pearson} correlation coefficient and the sample rank correlation, \chn{Spearman's rho}, for a two-dimensional vector with heavy-tailed marginals. In Section~\ref{sec:linear_model} we present a simple model with an explicit linear dependence and show that, when the sample size grows to infinity, then Pearson's correlation coefficient does not converge to a constant but rather to a random variable involving stable distributions. We also verify analytically and numerically that the rank correlation provides a consistent statistical estimator for this model. Next, in Section~\ref{sec:nonnegative} we prove that if random variables are heavy-tailed \chn{with infinite second moment} and non-negative, then the sample \chn{Pearson} correlation coefficient never converges to a negative value. Thus, such sequence will never be classified as `disassortative'.  This result is extended to \chn{sequences of graphs} in Section~\ref{sec-RGs}, where we also obtain quite general convergence criteria in the infinite volume limit for the Pearson's correlation coefficient and the Spearman's rho. In Section
\ref{sec-RG-exam} analytical results are provided for \chn{Pearson's} correlation coefficient and rank correlations in the configuration model and the Preferential Attachment model. \chn{We also present an adaptation of the configuration model that has strong negative degree-degree dependencies and prove that Spearman's rho converges to the theoretically justified negative value while Pearson's coefficient converges to zero. Furthermore, we construct an example, where Pearson's correlation coefficient converges to a random variable.} Numerical results are presented in Section~\ref{sec-num-res-RGs}. We close the paper in Section~\ref{sec:discussion} with a discussion on our results and possible extensions thereof.}

\section{Correlations between random variables}
\label{sec-corr}
In this section we introduce the dependency measures studied in this paper. \ch{We start with \chr{a general description of dependency measures for random vectors $(X,Y)$.} This will provide the necessary intuition and framework in order to understand what happens when $X$ and $Y$ are the degrees of neighboring nodes in a network graph.} We present \chr{Pearson's sample} correlation coefficient in Section \ref{sec:corrcoef}, and introduce Spearman's rho in Section \ref{sec:spearman}. \chn{In Section~\ref{sec:linear_model} we demonstrate an ill behaviour of \chrr{Pearson's} sample coefficient in a simple model with linear dependencies, and in Section~\ref{sec:nonnegative} we show that if $X$ and $Y$ are non-negative then the \chrr{Pearson's} sample \chrr{coefficient} cannot converge to a negative value.}

\subsection{Sample Pearson's correlation coefficient}
\label{sec:corrcoef}

The {Pearson correlation coefficient} $\rho$ for two random variables
$X$ and $Y$ \chr{with cumulative distribution functions $F_{\sss X}(\cdot)$ and $F_{\sss X}(\cdot)$, joint cumulative distribution function $F_{\sss X,Y}(\cdot,\cdot)$},  and $\Var(X),\Var(Y)<\infty$ is defined by
    	\eqn{
	\label{sample-corr-def}
	\rho=\frac{\expec[XY]-\expec[X]\expec[Y]}{\sqrt{\Var(X)}\sqrt{\Var(Y)}}.
    	}
By Cauchy-Schwarz, $\rho\in [-1,1]$, and $\rho$ measures the \emph{linear
dependence} between the random variables $X$ and $Y$.
We can approximate $\rho$ from a sample by computing the
\emph{sample correlation coefficient}
    	\eqn{
    	\label{rhon-def}
    	\rho_n=\frac{\frac{1}{n-1}\sum_{i=1}^n (X_i-\bar{X}_n)(Y_i-\bar{Y}_n)}
    	{S_n(X)S_n(Y)},
   	 }
where
    	\eqn{
	\label{sample-stats}
	\bar{X}_n=\frac{1}{n}\sum_{i=1}^n X_i,
    	\qquad
    	\bar{Y}_n=\frac{1}{n}\sum_{i=1}^n Y_i
    	}
denote the \emph{sample averages} of $(X_i)_{i=1}^n$ and
$(Y_i)_{i=1}^n$, while
    	\eqn{
    	\label{Sn-def}
    	S_n^2(X)=\frac{1}{n-1}\sum_{i=1}^n (X_i-\bar{X}_n)^2,
    	\qquad
    	S_n^2(Y)=\frac{1}{n-1}\sum_{i=1}^n (Y_i-\bar{Y}_n)^2
    	}
denote the \emph{sample variances.}  For i.i.d.\ sequences of random vectors $((X_i,Y_i))_{i=1}^n$  under the assumption of finite-variance random variables, i.e.,
 $\Var(X),\Var(Y)<\infty$, it is well known that
the estimator $\rho_n$ of $\rho$ is \emph{consistent,}
i.e.,
	\eqn{
	\label{consist-rho}
	\rho_n\convp \rho,
	}
where $\convp$ denotes convergence in probability. In practice,
however, we tend not to know whether $\Var(X),\Var(Y)<\infty$,
since $S_n^2(X)<\infty$ and $S_n^2(Y)<\infty$ clearly 
hold for any sample, and, therefore, one might be tempted to
always use $\rho_n$. Furthermore, by \ch{the} Cauchy-Schwarz \ch{inequality}, 
$\rho_n\in[-1,1]$ for every $n\geq 1$, which is part of the problem, because, for any sample, a value in 
$[-1,1]$ is produced, and no alarm bells start rinkling when $\rho_n$ is used inappropriately. In this paper we investigate the case $\Var(X),\Var(Y)=\infty$, and show that the use of
$\rho_n$ in this case, and in particular in scale-free random graphs, is uninformative. \ch{For example, in case of negative correlations $\rho_n$ converges to zero when $n\to\infty$, which makes it impossible to compare \chn{the data} of different sizes. Moreover, if correlations are positive, $\rho_n$ may even converge to a random variable, thus it can produce very different numbers for two random structures of the same size created by the same mechanism. We provide such examples for linearly dependent random variables in Section~\ref{sec:linear_model} and for random graphs in Section~\ref{sec:random_assortativity}.}

\subsection{Rank correlations}
\label{sec:spearman}

For two-dimensional data $((X_i, Y_i))_{i=1}^n$, let $r_i^{X}$ and $r_i^{Y}$ be the rank of an observation $X_i$ and $Y_i$, respectively,  when the sample values $(X_i)_{i=1}^n$ and
$(Y_i)_{i=1}^n$ are arranged in a descending order. The idea of rank correlations is in evaluating statistical dependences on the data $((r_i^X,r_i^Y))_{i=1}^n$, rather than on the original data $((X_i, Y_i))_{i=1}^n$. 
Rank transformation is convenient, in particular because, for continuous random variables, the two marginals of the resulting vector $(r^X_i,r_i^Y)$ are realizations of identical uniform distributions, implying many nice mathematical properties.

The statistical correlation coefficient for the ranks is known as Spearman's rho~\cite{Spearman1904}:
	\begin{equation}
	\label{eq:spearman}
	\rho_n^{\rm rank}=\frac{\sum_{i=1}^n(r_i^X-(n+1)/2)(r_i^Y-(n+1)/2)}
	{\sqrt{\sum_{i=1}^n(r_i^X-(n+1)/2)^2\sum_i^n(r_i^Y-(n+1)/2)^2}}
	\ch{=\frac{\frac{1}{n}\sum_{i=1}^n r_i^Xr_i^Y-((n+1)/2)^2)}{\frac{1}{12}\,(n^2-1)}}.
	\end{equation}
The mathematical properties of Spearman's rho have been extensively investigated in the literature. \ch{It is well known that} if $((X_i, Y_i))_{i=1}^n$ consists of independent realizations of $(X,Y)$, and the joint distribution \chr{cumulative} function of $X$ and $Y$ is continuous, then $\rho_n^{\rm rank}$ \ch{converges to a number that can be interpreted as its population value, see~\cite[Chapter~9]{Kendall1975}, \cite{Borkowf2002}:
	 \eqn{
	\label{consist-Spearman's}
	\rho_n^{\rm rank}\convp \rho^{\rm rank}=12\E(F_{\sss X}(X)F_{\sss Y}(Y))-3.
	}}
\ch{For completeness, we give a brief explanation of this formula. Observe that $F_{\sss X}(X)$ is the random variable that takes the value $F_{\sss X}(x)$ when $X=x$. If $X$ is continuous then $F_{\sss X}(X)$ has a uniform distribution on $[0,1]$:
	\begin{equation}
	\label{eq:uniform}
	F_{\sss X}(x)=\prob(X\le x)=\prob(F_{\sss X}(X)\le F_{\sss X}(x)).
	\end{equation}
Now take $F_{\sss X}(x)=t$ to obtain $\prob(F_{\sss X}(X)\le t)=t$, where $t$ can take any value in $[0,1]$. We note that this derivation holds for {\it any} continuous random variable $X$. We will use this many times throughout the paper. In particular, it follows that $\E(F_{\sss X}(X))=\chr{\E(F_{\sss Y}(Y))}=1/2$. Next, note that 
$r_i^X/n$ is an empirical estimator of $1-F_{\sss X}(x_i)$, where $x_i$ is the realized value of $X_i$. Moreover, 
	\[
	\E((1-F_{\sss X}(X))(1-F_{\sss Y}(Y)))
	=1-\E(F_{\sss X}(X))-\E(F_{\sss Y}(Y))+\E(F_{\sss X}(X)F_{\sss Y}(Y))=\E(F_{\sss X}(X)F_{\sss Y}(Y)).
	\]
Hence, the right-hand side of (\ref{eq:spearman}) is a statistical estimator of the last expression in (\ref{consist-Spearman's}).}

\label{tie-breaking}
For discrete random variables, the situation is more delicate, as the same values \chn{of $X$ and $Y$} may occur more than once. \ch{We resolve the ties randomly, using uniformisation as suggested in \cite{Mesfioui2005Spearman}. Formally, we replace the ranks of $\chr{((X_i,Y_i)_{i=1}^n}$ by the ranks of the random variables 
	\[
	\chr{((X^*_i,Y^*_i))_{i=1}^n=((X_i+U_i, Y_i+\chn{U'_i}))_{i=1}^n,}
	\] 
where $\chr{((U_i,U'_i))_{i=1}^n}$ \chr{is a sequence of $2n$} i.i.d.\ uniform variables on $(0,1)$. The random variables $\chr{X^*_i}$ and $\chr{Y^*_i}$ \chr{now} \emph{are} continuous. We denote their cumulative distribution functions by $F_{\sss X}^*$ and $F_{\sss Y}^*$. Note that if $X$ takes non-negative integer values then $F_{\sss X}^*$ can be seen as a linear interpolation of the cumulative probability $\prob(X<x)$, $x=0,1,2,\ldots$ because $\prob(X=x)=\prob(X^*\in [x,x+1))$.

Since $(X^*,Y^*)$ has a continuous distribution, the convergence result in (\ref{consist-Spearman's}) remains valid. Moreover, \cite{Mesfioui2005Spearman} gives the formula for $\rho^{\rm rank}$ in a discrete case, and \cite[Proposition~3.1]{Mesfioui2005Spearman} states that if $X,Y=0,1,\ldots$, then $(X,Y)$ and $(X^*,Y^*)$ have the same population value $\rho^{\rm rank}$:
	\begin{equation}
	\label{eq:proposition3-1}
	\rho^{\rm rank}=12\E(F^*_{\sss X}(X^*)F^*_{\sss Y}(Y^*))-3.
	\end{equation}
The comparison of different ways for resolving ties, and their effect on the resulting computation is an interesting topic, which is outside the scope of this work.} We refer to \cite{Nevslehova2007rank_correlations} for a general treatment of rank correlations for non-continuous distributions.

\subsection{Linear dependencies}
\label{sec:linear_model}

It is well known that $\rho$ in general measures {\it linear} dependence between two 
random variables. \ch{Therefore, before analyzing the behavior of $\rho_n$ in networks, we wish to illustrate that $\rho_n$ fails to capture the linear dependence between $X$ and $Y$ when \chrr{the variances of $X$ and $Y$ are infinite, i.e.,} ${\rm Var}(X),{\rm Var}(Y)=\infty$, even in a very straightforward case when the linear relation \chr{between} $X$ and $Y$ is explicitly defined. With this goal in mind, below we analyze the behavior of $\rho_n$} in the following linear model:
    \chn{\begin{eqnarray}
    \label{eq:X} 
	X&=\alpha_1\xi_1+\cdots+\alpha_m\xi_m,
	\qquad
    Y=\beta_1\xi_1+\cdots+\beta_m\xi_m,
    \end{eqnarray}}
where $\xi_j$, $j=1,\ldots,m$, are independent  identically distributed
(i.i.d.) non-negative random variables with regularly varying tail, and tail exponent $\gamma$. By definition, the
non-negative random variable $\xi$ is {\it regularly varying} with index
$\gamma>0$, if 
	\eqn{
    	\label{tail-U}
    	\P(\xi>x)=L(x)x^{-\gamma},\qquad x\ge 0,
    	}
where $\chrr{x\mapsto L(x)}$ is a slowly varying function, that is, for $u>0$,
$L(ux)/L(x)\to 1$ as $x\to\infty$, for instance, $L(x)$ may be
equal to a constant or $\log(x)$. Note that the random variables 
$X$ and $Y$ have the same distribution when $(\beta_1,\ldots, \beta_m)$
is a permutation of $(\alpha_1,\ldots, \alpha_m)$. 

When we take an i.i.d.\ sample of random variables $((X_i,Y_i))_{i=1}^n$ of random variables with
the above linear dependence, then Spearman's rho is consistent by \eqref{consist-Spearman's},
with a variance that converges to zero as $1/n$. 
For the sample correlation coefficient, consistency follows from \eqref{consist-rho} 
in the case where $\Var(\xi_i)<\infty$, but not when the $\xi_i$'s have infinite variance \chr{as we show below in detail}.
Our main result in this section is the following theorem:

\begin{Theorem}[Weak convergence of the \chn{sample Pearson's coefficient}]
\label{thm-corr-coef}
Let $((X_i,Y_i))_{i=1}^n$ be i.i.d.\ copies of the random variables
$(X,Y)$ in \eqref{eq:X}, and where
$(\xi_j)_{j=1}^m$ are i.i.d.\ random variables satisfying \eqref{tail-U}
with $\gamma\in (0,2)$, so that $\Var(\xi_j)=\infty$.
Then,
    \eqn{
    \label{rhon-asym}
    \rho_n\convd \rho\equiv \frac{\sum_{j=1}^m \alpha_\ch{j}\beta_\ch{j} Z_j}{\sqrt{\sum_{j=1}^m \alpha_\ch{j}^2 Z_j}
    \sqrt{\sum_{j=1}^m \beta_\ch{j}^2 Z_j}},
    }
where $(Z_j)_{j=1}^m$ are i.i.d.\ random variables having stable distributions with parameter $\gamma/2\in (0,1)$,
and $\convd$ denotes convergence in distribution.
In particular, 
$\rho$ has a density on $[-1,1]$. \ch{This density is strictly positive on $(-1,1)$ when
there exist $k,l$ such that $\alpha_k\beta_k<0<\alpha_l\beta_l$. Furthermore, the density
is positive on $(a,1)$ when $\alpha_k\beta_k\geq 0$ for every $k$, and on $(-1,-a)$ when $\alpha_k\beta_k\leq 0$ for every $k$}, where
\eqn{
    \label{def-interval-a}
    a=\inf_{z_1, \ldots, z_m\chr{\in \R}} \frac{\sum_{j=1}^m |\alpha_j\beta_j| z_j}{\sqrt{\sum_{j=1}^m \alpha_j^2 z_j}
    {\sqrt{\sum_{j=1}^m \beta_j^2z_j}}}\in (0,1).
    }
\end{Theorem}

Theorem \ref{thm-corr-coef} states that the sample correlation coefficient converges in distribution 
to a proper random variable, contrary to Spearman's rank correlation which converges in probability to a constant. In particular, this implies that when we have two independent samples, the sample correlation 
coefficient will give two rather distinct values, while Spearman's rank correlation will give two similar values.
We prove Theorem \ref{thm-corr-coef} in the remainder of this section.
In its proof, we need the following technical result:

\begin{Lemma}[Asymptotics of sums in stable domain]
\label{lem-sums-stable}
Let $(\xi_{i,j})_{i=1,2,\ldots,n, j=1,2}$ be i.i.d. random variables satisfying \eqref{tail-U} for some $\gamma\in (0,2)$. Then
there exists a sequence $a_n$ with $a_n=n^{2/\gamma}\ell(n)$, where $n\mapsto \ell(n)$ is
slowly varying, such that
    \eqn{
    \label{stable-conv}
    \frac{1}{a_n}\sum_{i=1}^n \xi_{i,1}^2\convd Z_1,
	\qquad
	\frac{1}{a_n}\sum_{i=1}^n \xi_{i,1}\xi_{i,2}\convp 0,
    }
where $Z_1$ is stable with parameter $\gamma/2$ and $\convp$ denotes
convergence in probability.
\end{Lemma}

\proof \ch{Let $F(x)=\prob(\xi\le x)$ be the cumulative} distribution function of $\xi$. 
\ch{In order to \chr{prove} the first statement in \eqref{stable-conv} we only need to note} that
the \chr{cumulative} distribution function of $\xi^2$ equals $x\mapsto F(\sqrt{x})$, which, by
\eqref{tail-U}, \ch{implies that $\xi^2$ is regularly varying. Thus, the first statement in \eqref{stable-conv} is in fact \chr{the classical convergence of infinite variance random variables with slowly varying distribution 
\chr{functions} to stable laws} (see e.g.\ \cite{GneKol68}), where $Z_1$ is a stable $\gamma/2$
random variable.} In particular, denoting $[1-F](x)=1-F(x)$, $x\ge 0$, we can identify $a_n=[1-F]^{-1}(1/n^2)$~\chr{\cite{BinGolTeu89}}.
\chr{Since $x\mapsto [1-F](x)$ is regularly varying with index $\gamma$,
$[1-F]^{-1}(1/n)$ is regularly varying with index $1/\gamma$ \cite{BinGolTeu89},
so that $a_n=[1-F]^{-1}(1/n^2)$ is regularly varying with index $2/\gamma$.}
To prove the second part of \eqref{stable-conv}, we write
    \eqn{
    \lbeq{[1-F]-infvar-ub}
    1-F(x)=
    \prob(\xi>x)\leq c' x^{-\gamma'},\qquad x\ge 0,
    }
which is valid for any $\gamma'\in (1,\gamma)$ by \eqref{tail-U} and Potter's theorem. 
We next study the \chr{cumulative} distribution function of $\xi_1\xi_2$
which we denote by $H$, 
where $\xi_1$ and $\xi_2$ are two independent copies of the
random variable $\xi$. When $F$
satisfies \eqref{[1-F]-infvar-ub}, then it is not hard to see
that there exists a $C>0$ such that
    \eqn{
    \label{1-G-bd}
    1-H(u)\leq C (1+\log{u}) u^{-\gamma'}.
    }
Indeed, assume that $F$ has a density $f(w)=c w^{-(\gamma'+1)}$, for $w\geq 1$.
Then,
    $$
    1-H(u)
    =\int_1^{\infty} f(w) [1-F](u/w)dw.
    $$
Clearly, $1-F(w)= c' w^{-\gamma'}$ for $w\geq 1$ and $1-F(w)=1$ otherwise. Substitution of this
yields
    $$
    1-H(u)
    \le cc' \int_1^{u} w^{-(\gamma'+1)} (u/w)^{-\gamma'}dw+c\int_u^{\infty}
    w^{-(\gamma'+1)}\,dw
    \leq C (1+\log{u}) u^{-\gamma'}.
    $$
When $F$ satisfies \eqref{[1-F]-infvar-ub}, then $\xi_1$ and $\xi_2$
are stochastically upper bounded by \chn{$\hat{\xi}_1$ and $\hat{\xi}_2$} with
\chr{cumulative} distribution function $\hat{F}$ satisfying $1-\hat{F}(w)=c'w^{-\gamma'}\vee 1$, 
where $(x\vee y)=\max\{x,y\}$, and the claim in \eqref{1-G-bd} follows from the above computation.

By the bound in \eqref{1-G-bd}, the random variables $\xi_{i,1}\xi_{i,2}$
are stochastically bounded from above by random variables $P_i$
that are in the domain of attraction of a stable $\gamma'$
random variable. As a result, there exists $b_n=n^{1/\gamma'}\ell'(n)$, where $n\mapsto \ell'(n)$ is
slowly varying, such that
    \[
    \frac{1}{b_n}\sum_{i=1}^n P_i\convd W,
    \]
where $W$ is stable $\gamma'$. {By choosing $\gamma'>\gamma/2$}, we get $b_n/a_n\to 0$, so we obtain
the second statement in \eqref{stable-conv}.
\qed
\bigskip

\proof[Proof of Theorem~\ref{thm-corr-coef}]  We start by noting that
    \eqn{
    \label{rhon-def-rep}
    \rho_n=\frac{\frac{1}{n-1}\sum_{i=1}^n (X_iY_i-\bar{X}_n\bar{Y}_n)}
    {S_n(X)S_n(Y)},
    }
and
    \eqn{
    \label{Sn-def-rep}
    S_n^2(X)=\frac{1}{n-1}\sum_{i=1}^n(X_i^2-\bar{X}_n^2),
    \qquad
    S_n^2(Y)=\frac{1}{n-1}\sum_{i=1}^n (Y_i^2-\bar{Y}_n^2).
    }
We continue to identify the asymptotic behavior of
    \[
    \sum_{i=1}^nX_i^2,
    \qquad
    \sum_{i=1}^nY_i^2,
    \qquad
    \sum_{i=1}^nX_iY_i.
    \]
Let $[n]$ denote the set of integers $\{1,2,\ldots,n\}$. The distribution of  $((X_i,Y_i))_{i=1}^n$ is described in
terms of an array $(\xi_{i,j})_{i\in [n], j\in [m]}$, which are i.i.d.\
copies of a random variable $\xi$. In terms of these random variables, we can identify
    \eqn{
    \label{sum-prods}
    \sum_{i=1}^nX_iY_i
    =\sum_{j=1}^m \alpha_j\beta_j \Big(\sum_{i=1}^n \xi_{i,j}^2\Big)
    +\sum_{j_1\neq j_2=1}^m \alpha_{{j_1}}\beta_{{j_2}} \Big(\sum_{i=1}^n \xi_{i,j_1}\chrr{\xi}_{i,j_2}\Big).
    }
The sums $\sum_{i=1}^n \xi_{i,j}^2$ are i.i.d.\ for different $j\in \{1, \ldots, m\}$,
and by Lemma~\ref{lem-sums-stable}, $\sum_{i=1}^n \xi_{i,j_1}\xi_{i,j_2}$ is of a smaller
order. Hence, from \eqref{sum-prods} we obtain that
    \eqn{
    \frac{1}{a_n} \sum_{i=1}^nX_iY_i
    \convd \sum_{j=1}^m \alpha_j\beta_j Z_j.
    }
Therefore, by
taking $\alpha=\beta$, we also obtain
    \eqn{
    \label{conv-sec-mom}
    \frac{1}{a_n} \sum_{i=1}^nX_i^2
    \convd \sum_{j=1}^m \alpha_j^2 Z_j,
    \qquad
    \frac{1}{a_n} \sum_{i=1}^nY_i^2
    \convd \sum_{j=1}^m \beta_j^2 Z_j,
    }
and \chrr{the convergence holds} \emph{simultaneously}. As a result, \eqref{rhon-asym}
follows. It remains to establish the properties of the limiting random variable
$\rho$ in \eqref{rhon-asym}. 

\chn{The density of $Z_i$ is strictly positive on $(0,\infty)$. Note that rescaling $z_j=cz_j$  $j=1,\ldots,m$, in (\ref{def-interval-a}), does not change the value of $a$. In particular, we can choose $c=(\max\{z_1,z_2,\ldots,z_m\})^{-1}$. If there exist $k$ and $l$ such that
$\alpha_k\beta_k<0<\alpha_l\beta_l$ then 
the density of $\rho$ is strictly positive on $(-1,1)$. Indeed, with positive probability $\rho$ can be arbitrarily close to $-1$ if $Z_k=\max\{Z_1,\ldots,Z_m\}$ and $Z_j/Z_k$, $j\ne k$ are sufficiently small. Similarly, if $Z_l=\max\{Z_1,\ldots,Z_m\}$ then with positive probability, $\rho$ can be arbitrarily close to $1$. \ch{Now assume that $\alpha_k\beta_k\geq 0$ for every $k$. In this case}, the density of $\rho$ is strictly positive on the support of $\rho$, which is $(a,1)$, with \chn{$a$ as in (\ref{def-interval-a})}. Analogously, when $\alpha_k\beta_k\leq 0$ then $\rho$ cannot be positive, and has a density on $(-1,-a)$.    }
\qed

\bigskip

\noindent
\chn{{\bf Numerical example.} In order} to illustrate the result of Theorem \ref{thm-corr-coef}, consider the example 
with $\xi_j$'s 
from a Pareto distribution satisfying $\P(\xi>x)=1/x^{1.1}$, $x\ge 1$, so $L(x)=1$ and $\gamma=1.1$  
in \eqref{tail-U}. The exponent $\gamma=1.1$ is as observed for the World Wide Web~\cite{Broder00}. In \eqref{eq:X}, we choose $m=3$ and 
$\alpha_i$, $\beta_i$, $i=1,2,3$, as specified in Table~\ref{tab:linear_model}. 
We generate $N$ data samples $((X_i,Y_i))_{i=1}^n$ and compute $\rho_n$ and $\rho_n^{\rm rank}$ for each of the $N$ samples.
Thus, we obtain the vectors $(\rho_{n,j})_{j=1}^N$ and $(\rho^{\rm rank}_{n,j})_{j=1}^N$ of $N$ independent realizations for $\rho_n$ and $\rho_n^{\rm rank}$, respectively, where the sub-index $j=1,\ldots,N$ denotes the $j$th realization of $((X_i,Y_i))_{i=1}^n$. We then compute
	\begin{align}
	\label{eq:EN}
	\E_N(\rho_n)=\frac{1}{N}\sum_{j=1}^N \rho_{n,j},
    	\quad &\E_N(\rho^{\rm rank}_n)=\frac{1}{N}\sum_{j=1}^N \rho^{\rm rank}_{n,j};\\
	\label{eq:sigmaN}
	\sigma_N(\rho_n)=\sqrt{\frac{1}{N-1}\sum_{j=1}^N (\rho_{n,j}-\E_N(\rho_n))^2},\quad
	&\sigma_N(\rho^{\rm rank}_n)=\sqrt{\frac{1}{N-1}\sum_{j=1}^N (\rho^{\rm rank}_{n,j}-
	\E_N(\rho^{\rm rank}_n))^2}.
	\end{align}

\begin{figure}[ht]
\centerline{\includegraphics[height=1.2in]{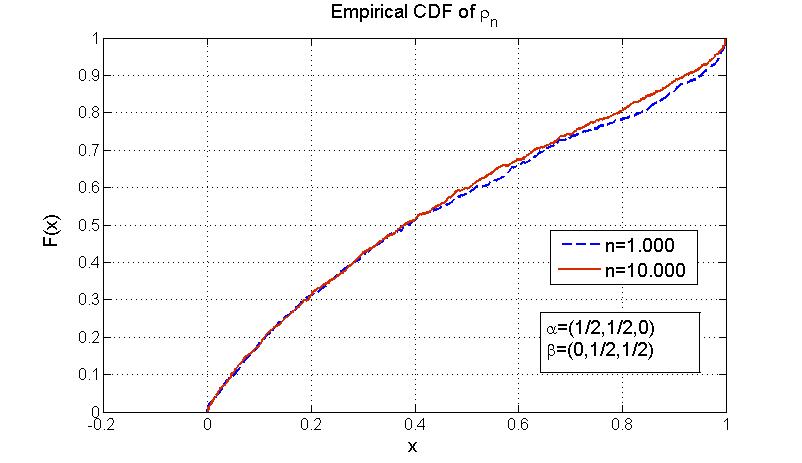}\hspace{.1in}
\includegraphics[height=1.2in]{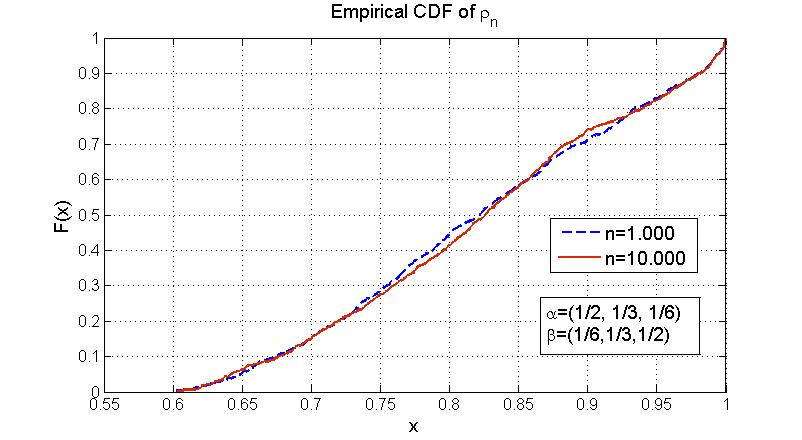}\hspace{.1in}
\includegraphics[height=1.2in]{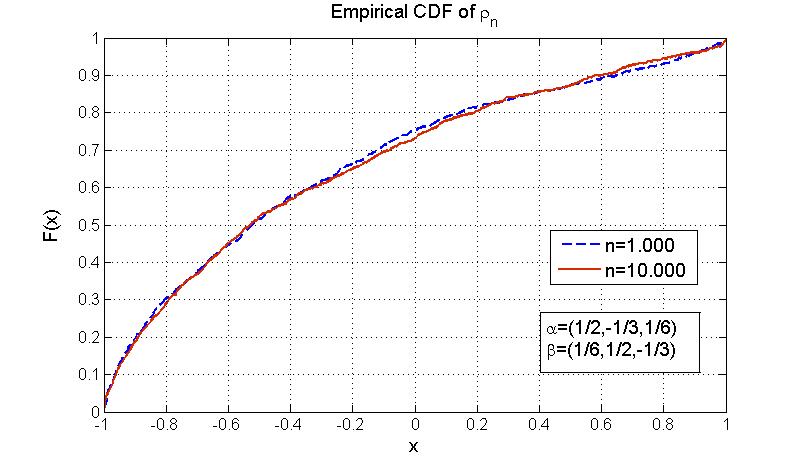}
}
\caption{\small The empirical distribution function $F_N(x)=\P(\rho_n\le x)$ for the $N=1.000$ observed values of $\rho_n$ ($n=1.000$, $n=10.000$), in the case of linear dependence \eqref{eq:X}.}
\label{fig:rho}
\end{figure}

The results are presented in Table~\ref{tab:linear_model}. We clearly see that $\rho_n$ has a significant standard deviation, of which estimators are similar for different values of $n$. This means that in the limit as $n\to\infty$, $\rho_n$ is a random variable with a significant spread in its values, as stated in Theorem~\ref{thm-corr-coef}. Thus, by evaluating $\rho_n$ for one sample $((X_i,Y_i))_{i=1}^n$ we will obtain a random number, even when $n$ is huge. 
The convergence to a non-trivial distribution is directly 
seen in Figure~\ref{fig:rho} because the plots for the two values of $n$ almost coincide. Note 
that in all cases, the density is fairly uniform, ensuring a comparable probability 
for all feasible values and rendering the value obtained in a specific realization
even more uninformative.  
\begin{table}[htb]%
{\small\centerline{
\begin{tabular}{|c|l|c|c|c|c|}
\hline
&$N$&$10^3$&\multicolumn{3}{|c|}{$10^2$}\\
\cline{2-6}
Model parameters&$n$&$10^2$&$10^3$&$10^4$&$10^5$\\
\hline
&$\E_N(\rho_n)$&0.4395&0.4365&0.4458&0.4067\\
$\alpha=(1/2,1/2,0)$&$\sigma_N(\rho_n)$&0.3399&0.3143&0.3175&0.3106\\
\cline{2-6}
$\beta=(0,1/2,1/2)$&$\E_N(\rho^{\rm rank}_n)$&0.4508&0.4485&0.4504&0.4519\\
&$\sigma_N(\rho^{\rm rank}_n)$&0.0922&0.0293&0.0091&0.0033\\
\hline
&$\E_N(\rho_n)$&0.8251&0.7986&0.8289&0.8070\\
$\alpha=(1/2,1/3,1/6)$&$\sigma_N(\rho_n)$&0.1151&0.1125&0.1108&0.1130\\
\cline{2-6}
$\beta=(1/6,1/3,1/2)$&$\E_N(\rho^{\rm rank}_n)$&0.8800&0.8850&0.8858&0.8856\\
&$\sigma_N(\rho^{\rm rank}_n)$&0.0248&0.0073&0.0023&0.0007\\
\hline
&$\E_N(\rho_n)$&-0.3052&-0.3386&-0.3670&-0.3203\\
$\alpha=(1/2,-1/3,1/6)$&$\sigma_N(\rho_n)$&0.6087&0.5841&0.5592&0.5785\\
\cline{2-6}
$\beta=(1/6,1/2,-1/3)$&$\E_N(\rho^{\rm rank}_n)$&-0.3448&-0.3513&-0.3503&-0.3517\\
&$\sigma_N(\rho^{\rm rank}_n)$&0.1202&0.0393&0.0120&0.0034\\
\hline
\end{tabular}}
\caption{\small Estimated mean and standard deviation of $\rho_n$ and $\rho^{\rm rank}_n$ in $N$ samples with linear dependence (\ref{eq:X}), $\P(\xi>x)=x^{-1.1}$, $x\ge 1$.}
\label{tab:linear_model}
}
\end{table}

On the other hand, from Table~\ref{tab:linear_model} we clearly see that the behaviour of the rank correlation is exactly as we can expect from a good statistical estimator.  The obtained average values are consistent while the standard deviation of $\rho_n^{\rm rank}$ decreases approximately as $1/\sqrt{n}$ as $n$ grows large. Therefore, $\rho_n^{\rm rank}$ converges to a deterministic number.

\subsection{Sample \chn{Pearson's correlation coefficient} for non-negative variables}
\label{sec:nonnegative}
We proceed by investigating correlations between non-negative heavy-tailed random variables. Our main result in this section shows that the correlation coefficient is asymptotically non-negative:

\begin{Theorem}[Asymptotic non-negativity of \chn{the sample Pearson's coefficient} for positive r.v.'s]
\label{thm-corr-coef-pos}
Let $((X_i,Y_i))_{i=1}^n$ be i.i.d.\ copies of non-negative random variables
$(X,Y)$, where $X$ and $Y$ satisfy
    \eqn{
    \label{tail-XY}
    \P(X>x)=L_{\sss X}(x)x^{-\gamma_{\sss X}},
    \qquad
    \P(Y>y)=L_{\sss Y} (y)y^{-\gamma_{\sss Y}},\qquad x,y\ge 0,
    }
with $\gamma_{\sss X}, \gamma_{\sss Y}\in (0,2)$, so that $\Var(X)=\Var(Y)=\infty$.
Then, any limit point of the sample \chrr{Pearson} correlation coefficient is non-negative.
\end{Theorem}

We illustrate Theorem \ref{thm-corr-coef-pos} with a useful example. Let $(\xi_i)_{i=1}^n$
be a sequence of i.i.d.\ random variables satisfying \eqref{tail-U} for some $\gamma\in (0,2)$,
and where $\xi\geq 0$ a.s. Let $(X,Y)=(0,2\xi)$ with probability $1/2$ and
$(X,Y)=(2\xi,0)$ with probability $1/2$. Then, $XY=0$ a.s., while $\expec[X]=\expec[Y]=\expec[\xi]$
and $\Var(X)=\Var(Y)=2\expec[\xi^2]-\expec[\xi]^2=2\Var(\xi)+\expec[\xi]^2$.
\chr{By Theorem \ref{thm-corr-coef-pos}, $\rho_n\convp0$ when $(\xi_i)_{i=1}^n$
is a sequence of i.i.d.\ non-negative random variables satisfying \eqref{tail-U} for some $\gamma\in (0,2)$, which is not appropriate as $(X,Y)$ are highly negatively dependent. When $\gamma>2$, this anomaly 
does not arise, since, if $\Var(\xi)<\infty$,}
    \eqn{
    \label{rhon-pos-fin-var}
    \rho_n\convp \rho=-\frac{\expec[\xi]^2}{2\Var(\xi)+\expec[\xi]^2}
    \in (-1,0).
    }
The asymptotics in \eqref{rhon-pos-fin-var} are quite reasonable, since the random
variables $(X,Y)$ are highly negatively dependent: When $X>0$, $Y$ must be equal to 0,
and vice versa. 

Table~\ref{tab:non-neg} shows the empirical mean and standard deviation of the estimators $\rho_n$ and $\rho_n^{\rm rank}$. Here $\P(\xi>x)=x^{-1.1}$, $x\ge 1$, as in Table~\ref{tab:linear_model}.
As predicted by Theorem~\ref{thm-corr-coef-pos}, the sample correlation coefficient (assortativity) converges to zero as $n$ grows large, while $\rho_n^{\rm rank}$ consistently shows a clear negative dependence, and the precision of the estimator improves as $n\to\infty$.
This explains why strong disassortativity is not observed in large samples of \chrr{non-negative} power-law data.

\begin{table}%
{
{\small\centerline{
\begin{tabular}{|l|c|c|c|c|c|}
\hline
$N$&\multicolumn{2}{|c|}{$10^3$}&\multicolumn{3}{|c|}{$10^2$}\\
\hline
$n$&$10$&$10^2$&$10^3$&$10^4$&$10^5$\\
\hline
$\E_N(\rho_n)$&-0.4833&-0.1363&-0.0342&-0.0077&-0.0015\\
$\sigma_N(\rho_n)$&0.1762&0.0821&0.0245&0.0064&0.0011\\
\hline
$\E_N(\rho^{\rm rank}_n)$&-0.6814&-0.4508&-0.4485&-0.4504&-0.4519\\
$\sigma_N(\rho^{\rm rank}_n)$&0.1580&0.0283&0.0082&0.0024&0.0007\\
\hline
\end{tabular}}}
\caption{\small The  mean and standard deviation of $\rho_n$ and $\rho_n^{\rm rank}$ in $N$ simulations of $((X_i,Y_i))_{i=1}^n$, where $X=2\chrr{\xi}I$, $Y=2\chrr{\xi}(1-I)$, $I$ is a Bernoulli$(1/2)$ random variable, $\P(\xi>x)=x^{-1.1}$, $x\ge 1$.}
\label{tab:non-neg}}
\end{table}
\medskip

\noindent
{We next prove Theorem \ref{thm-corr-coef-pos}:\\
{\it Proof of Theorem \ref{thm-corr-coef-pos}.}
Clearly $\sum_{i=1}^nX_iY_i\geq 0$ when $X_i\geq 0, Y_i\geq 0$, so that
    \[
    \rho_n\geq -\frac{\frac{1}{n-1}\sum_{i=1}^n\bar{X}_n\bar{Y}_n}
    {S_n(X)S_n(Y)}
    =-\frac{n}{n-1} \frac{\bar{X}_n}{S_n(X)}\frac{\bar{Y}_n}{S_n(Y)}.
    \]
{It remains to show that if $\Var(X)=\infty$, then $\bar{X}_n/S_n(X)\convp 0$. Indeed, if $\gamma\in (1,2)$ then $\bar{X}_n\convp \expec[X]<\infty$ by the strong
law of large numbers. When $\gamma\in(0,1]$, instead, then $X$ is in the domain of attraction of a $\gamma$ stable random variable, hence $\bar{X}_n$, loosely speaking, it scales as $n^{1/\gamma_{\sss X}-1}$. Further, from \eqref{tail-XY} and Lemma~\ref{lem-sums-stable} it follows that $S_n(X)$ scales as $n^{2/\gamma_{\sss X}-1}$, in particular, $\bar{X}_n/S_n(X)\convp 0$ for all $\gamma\in \chn{(0,2)}$. }
\qed

\section{Applications to networks}
\label{sec-RGs}
 \chn{In real-world networks \chrr{it} is particularly important to measure degree-degree dependencies for neighboring \chrr{vertices.}} We refer to \cite{Newman} for an extensive introduction to networks, 
their empirical properties and models for them. In Section \ref{sec:definitions} \chrr{below,} we start with \chrr{the} formal definition of \chn{Pearson's correlation coefficient (which was termed \chrr{the} {\it assortativity coefficient} in \cite{Newman2002assortative}),} and Spearman's rho in the network context.
Next, in Section~\ref{sec-no-disass-RG-seq}
we show that all limit points of \chn{Pearson's} coefficients for sequences of growing scale-free
random graphs with power-law exponent $\gamma<3$ are non-negative,
a result that is similar in spirit to Theorem \ref{thm-corr-coef-pos}.
In Section \ref{sec-conv-rank-assor-gen}, we state general convergence conditions
for both \chn{Pearson's} correlation coefficient as well as Spearman's rho. 
}

\subsection{Definitions and notations}
\label{sec:definitions}

We start by introducing some notation. \chn{Let $G=(V,E)$ be an undirected random graph. For a {\it directed} edge $e=(u,v)$,
we write $\underline{e}=u, \overline{e}=v$ and we denote the set of directed edges in $E$ \chrr{by $E'$} (so that $|E'|=2|E|$), and $D_v$ is the degree of vertex $v\in V$. In general, $D_v$ is a random variable.}

The assortativity coefficient of $G$ is equal to
(see, e.g., \cite[(4)]{Newman2002assortative})
    \begin{equation}
    \label{eq:def_pearson}
    \rho(G)=\frac{\frac{1}{|E'|}\sum_{(u,v)\in E'} D_uD_v-\Big(\frac{1}{|E'|} \sum_{(u,v)\in E'} \tfrac12 (D_u+D_v)\Big)^2}
    {\frac{1}{|E'|}\sum_{(u,v)\in E'} \tfrac12(D_u^2+D_v^2) -\Big(\frac{1}{|E'|} \sum_{(u,v)\in E'} \tfrac12 (D_u+D_v)\Big)^2}.
    \end{equation}
Note that the assortativity coefficient in (\ref{eq:def_pearson}) is equal to the sample correlation coefficient, where $((D_u, D_v))_{(u,v)\in E'}$ represent a sequence of non-negative random variables, as studied in Theorem~\ref{thm-corr-coef-pos}. However, $((D_u, D_v))_{(u,v)\in E'}$ are not independent, so that we may not immediately apply the previous theory. \chrr{Theorem~\ref{thm-corr-coef-pos-scale-free} below is the analogue of Theorem~\ref{thm-corr-coef-pos} in the network context, and we give a formal proof of it below.}

Let us now introduce Spearman's rho in $G$ that we denote by $\rho^{\rm rank}(G)$. In accordance to the original definition of Spearman's rho, $\rho^{\rm rank}(G)$ is the correlation coefficient of 
the sequence of random variables $(R_{\underline{e}},R_{\overline{e}})$,
where $e$ is a uniformly chosen directed edge $(u,v)$ from $E_n'$. We let $R_{\underline{e}}$ and $R_{\overline{e}}$ be the \emph{rank} of \chn{respectively 
$D_{\underline{e}}+U_e$ and $D_{\overline{e}}+U_e'$ in the sequences $(D_{\underline{e}}+U_{e})_{e\in E'_n}$ and $(D_{\overline{e}}+U_{e}')_{e\in E_n'}$. 
Here, as discussed on page \pageref{tie-breaking}, \chr{$(U_{e})_{e\in E_n'}$ and $(U_{e}')_{e\in E_n'}$}
are i.i.d.\ sequences of uniform $(0,1)$ random variables.
\chrr{Then,} \chr{Spearman's} rank correlation coefficient is defined as follows:}
\chn{\eqn{
    	\label{eq:rhoG_rank}
    	\rho^{\rm rank}(G)
	=\frac{\frac{1}{|E'|}\sum_{e\in E'} R_{\underline{e}}R_{\overline{e}}-(|E'|+1)^2/4}{(|E'|^2-1)/12}.}
	}

\subsection{No disassortative scale-free random graph sequences}
\label{sec-no-disass-RG-seq}
{We compute that 
    \eqn{
	\label{relation-moment-degrees}
    \frac{1}{|E'|} \sum_{(u,v)\in E'} \tfrac12 (D_u+D_v)
    =\frac{1}{|E'|} \sum_{v\in V} D_v^2,
    \qquad
    \frac{1}{|E'|}\sum_{(u,v)\in E'} \tfrac12(D_u^2+D_v^2)
    =\frac{1}{|E'|} \sum_{v\in V} D_v^3.
    }}
Thus, $\rho(G)$ can be written as
    \eqn{
    \label{eq:rhoG}
    \rho(G)=\frac{\sum_{(u,v)\in E'} D_uD_v-\frac{1}{|E'|} \Big(\sum_{v\in V} D_v^2\Big)^2}
    {\sum_{v\in V} D_v^3 -\frac{1}{|E'|} \Big(\sum_{v\in V} D_v^2\Big)^2}.
    } 

Consider a sequence of graphs $(G_n)_{n\geq 1}$, where $G_n=(V_n, E_n)$ and
$n$ denotes the number of vertices $n=|V_n|$ in the graph. Since many real-world networks
are quite large, we are interested in the behavior of
$\rho(G_n)$ as $n\rightarrow \infty$. Note that this
discussion applies both to sequences of real-world networks 
of increasing size, as well as to graph sequences of random 
graphs. We start by generalizing Theorem \ref{thm-corr-coef-pos} 
to this setting:

\begin{Theorem}[Asymptotic \chn{non-negativity of Pearson's coefficient} in scale-free graphs]
\label{thm-corr-coef-pos-scale-free}
Let $(G_n)_{n\geq 1}$  be a sequence of graphs of size $n$
satisfying that there exist $\gamma\in (1,3)$ and $0<c<C<\infty$ such that
$cn\leq |E|\leq Cn$,  $c n^{1/\gamma}\le \max_{v\in V_n} D_v \leq C n^{1/\gamma}$ and $cn^{(2/\gamma)\vee 1}\leq \sum_{v\in V_n} D_v^2 \leq Cn^{(2/\gamma)\vee 1}$.  
Then, any limit point of \chn{Pearson's} correlation coefficient
$\rho(G_n)$ is non-negative.
\end{Theorem}
\medskip

\chr{In the next section, we give several examples where Theorem \ref{thm-corr-coef-pos-scale-free}
applies and yields results that are not sensible. The \chn{powerful feature of} Theorem \ref{thm-corr-coef-pos-scale-free} is that it applies to \emph{all} graphs, not just realizations of certain random graphs.}

\proof
We note that $D_v\geq 0$ for every $v\in V$,
so that, from (\ref{eq:rhoG})
    \eqn{\label{eq:rho-}
    \rho(G_n)\geq {\rho}^-(G_n)\equiv 
    -\frac{\frac{1}{|E'|}\Big(\sum_{v\in V} D_v^2\Big)^2}
    {\sum_{v\in V} D_v^3-\frac{1}{|E'|}\Big(\sum_{v\in V} D_v^2\Big)^2}.
    }
By assumption, $\sum_{v\in V} D_v^3\geq (\max_{v\in [n]} D_v)^3
\geq c^3 n^{3/\gamma}$, whereas $\frac{1}{|E'|}\Big(\sum_{v\in V} D_v^2\Big)^2
\leq (C^2/c) n^{2(2/\gamma\vee 1)-1}=(C^2/c) n^{[(4/\gamma-1)\vee 1]}$. Since $\gamma\in (1,3)$ we have $(4/\gamma-1)\vee 1<3/\gamma$, so that
	\[
	\frac{\sum_{v\in V} D_v^3}{\frac{1}{|E'|}\Big(\sum_{v\in V} D_v^2\Big)^2}
	\rightarrow \infty.
	\]
\ch{Hence, ${\rho}^-(G_n)\to 0$ as $n\to\infty$.} This proves the claim.
\qed
\medskip

In the literature, many examples are reported of real-world networks
where the degree distribution \ch{closely follows} a power law with $\gamma$ in $(1,3)$, see e.g., \cite[Table I]{Albert2002stat_mech}
or \cite[Table I]{Newman03survey}. \chn{Let $D$ be such \chrr{a} power-law random variable, and \chrr{denote} $\mu_p=\expec[D^p]$ for $p\in(0,\gamma)$.} In that case one 
can expect that
    \[
    |E'|= \sum_{v\in V} D_v\sim \mu_1 n,
    \]
while $\max_{v\in V}D_v\sim n^{1/\gamma}$, and
    \eqn{
	\label{power-degree-conv}
    \frac{1}{n} \sum_{v\in V} D_v^p
    \sim 
    \begin{cases}
    \mu_p &\text{when }\gamma>p,\\
    \chn{C_p} n^{p/\gamma-1} &\text{when }\gamma<p.
    \end{cases}
    }  \ch{Of course, the convergence in \eqref{power-degree-conv} depends sensitively on the occurrence of 
large degrees. However, intuitively it can be explained as follows. When}\chr{ 
	\[
	\frac{1}{n} \sum_{v\in V} \indic{D_v\geq k} =C'  k^{-\gamma}(1+o(1))
	\]
for all $k$ for which $k^{-\gamma}\gg 1/n$ so that $k\ll n^{1/\gamma}$, then
	\[
	\frac{1}{n} \sum_{v\in V} D_v^p=\sum_{k\geq 1} (k^p-(k-1)^p) \frac{1}{n} \sum_{v\in V}\indic{D_v\geq k}
	\chn{\approx C''  \sum_{k=1}^{n^{1/\gamma}} k^{p-1-\gamma}= C_p  n^{p/\gamma-1},}
	\]}
\chn{where $C''$ and $C_p$ are appropriately chosen constants.}
In particular, the conditions of Theorem \ref{thm-corr-coef-pos-scale-free} hold and 
${\rho}^-(G_n)\rightarrow 0 \quad \text{when }\gamma<3.$
Thus, the asymptotic degree-degree correlation of the graph sequence      
$(G_n)_{n\geq 1}$ is non-negative. As a result, when the power-law exponent satisfies $\gamma<3$
there exist no  scale-free graph 
sequences \chn{that will be identified as disassortative by Pearson's coefficient}.
We next investigate a general theorem that allows us 
to identify the limit of Spearman's rho and \chn{Pearson's} coefficient
for many random graph models.

\subsection{Convergence conditions for \chn{degree-degree dependency measures}}
\label{sec-conv-rank-assor-gen}
Let $(G_n)_{n\geq 1}$  be \chn{again} a sequence of graphs of size $n$,
where $G_n=(V_n,E_n)$, $|V_n|=n$. We write $\expec_n$ for the conditional expectation given the graph $G_n$
(which in itself is random, so that we are \emph{not} taking the
expectation w.r.t.\ $G_n$). Consider a random vector $(X,Y)=(D_{\underline{e}},D_{\overline{e}})$ where $e$ is chosen uniformly at random from $E'$. \ch{Recall that for a discrete random variable
$X$, $F_{\sss X}$ denotes its \ch{cumulative} distribution function, and 
$F^*_{\sss X}$ denotes the \ch{cumulative} distribution function of $X^*=X+U$, where $U$ is an independent 
uniform random variable on $(0,1)$. Then $F^*_{\sss X}(X^*)$ has a 
uniform distribution on $(0,1)$, see (\ref{eq:uniform}).} Our main result to identify the limits
of Spearman's rho as given by (\ref{eq:rhoG_rank}) and Pearson's coefficient is the following theorem:

\begin{Theorem}[Convergence criteria for \chn{degree-degree dependency measures}]
\label{thm-conv-rho-assor}
Let $(G_n)_{n\geq 1}$  be a sequence of random graphs of size $n$,
where $G_n=(V_n,E_n)$, $|V_n|=n$.  
Let $(X_n,Y_n)$ be the degrees on both sides of a uniform directed edge $e\in E'_n$.
Suppose that for every bounded continuous $h\colon \R^2\to \R$,
	\eqn{
	\label{conv-prob-degree}
	\expec_n[h(X_n,Y_n)]\convp \expec[h(X,Y)],
	}
where the r.h.s.\ is non-random. Then \\
(a)
	\eqn{\label{eq:spearman_convergence}
	\rho^{\rm rank}(G_n)\convp 12\E(F^*_{\sss X}(X^*)F^*_{\sss X}(Y^*))-3=\rho^{\rm rank},
	}
where \ch{$X^*=X+U$, $Y^*=Y+U'$,} $U$ and $U'$ are independent random variables on $(0,1)$, also independent of $X$ and $Y$, and $F^*_{\sss X}(\cdot)$ is the cumulative distribution function of $X^*$;\\
(b) \chr{when} we further suppose that $\expec_n[X_n^2]\convp\expec[X^2]<\infty$, and $\Var(X)>0$,
\chr{then \chrr{also}}
	\eqn{
	\rho(G_n)\convp \rho=\frac{\Cov(X,Y)}{\Var(X)}.
	}
\end{Theorem}

We remark that when $G_n$ is a random graph, then 
$\rho^{\rm rank}(G_n)$ and $\rho(G_n)$ are random variables. 
Equation \eqref{conv-prob-degree}
implies that the \emph{distribution} of the degrees on either side of an edge 
converges in probability to a deterministic limit, which can be interpreted
as the statement that the degree distribution converges to a deterministic limit.
The limits of $\rho^{\rm rank}(G_n)$ and $\rho(G_n)$ only depend on the
limiting degree distribution, where $\rho^{\rm rank}(G_n)$ \emph{always}
converges, while $\rho(G_n)$ can only be proved to converge when 
its limit \chn{is well defined}. We further note that \eqref{conv-prob-degree} 
is equivalent to showing that 
	\eqn{
	\label{conv-prob-degree-rep}
	\#\{e=(u,v)\in E'_n\colon (D_u, D_v)=(k,l)\}/|E_n'|\convp \prob(X=k,Y=l).
	} 
Condition \eqref{conv-prob-degree-rep} will be simpler to 
verify in practice. We emphasize that we study \emph{undirected} graphs but we work with \emph{directed} edges $e=(u,v)$,
which we vary over the whole set of edges, in such a way that $(u,v)$
and $(v,u)$ contribute as different edges. In particular, 
the marginal distributions of $X_n$ and $Y_n$
and consequently of $X$ and $Y$, are the same.
We next prove Theorem \ref{thm-conv-rho-assor}:

\proof We start with part (a). 
The sequence $(R_{\underline{e}}/|E_n'|,R_{\overline{e}}/|E_n'|)$ is a bounded
sequence of two-dimensional random variables. Let $F_{n,{\sss X}}$
denote the empirical \ch{cumulative} distribution \chr{function} of $(D_{\underline{e}})_{e\in E_n'}$ 
(which equals that of $(D_{\overline{e}})_{e\in {E_n'}}$), and let
\ch{$F^*_{n,{\sss X}}$} denote the empirical \ch{cumulative} distribution functions of 
$(D_{\underline{e}}+U_{e})_{e\in E'_n}$ (which equals that of
$(D_{\overline{e}}+U'_{e})_{e\in E'_n}$), where $(U_e)_{e\in E'_n}$, $(U'_e)_{e\in E'_n}$ are independent sequences of i.i.d uniform $(0,1)$ random variables.  Then, we can rewrite, with $\ell_n=|E_n'|$,
	\eqn{
	(R_{\underline{e}},R_{\overline{e}})
	=
	\big((\lceil \ell_n F^*_{n,{\sss X}}(D_{\underline{e}}+U_{e})\rceil, 
	\lceil \ell_n F^*_{n,{\sss X}}(D_{\overline{e}}+U'_{e})\rceil\big).
	}
In particular,
	\eqn{
	(R_{\underline{e}}/\ell_n,R_{\overline{e}}/\ell_n)
	=
	\big(\lceil \ell_n F^*_{n,{\sss X}}(D_{\underline{e}}+U_{e})\rceil/\ell_n, 
	\lceil \ell_n F^*_{n,{\sss X}}(D_{\overline{e}}+U'_{e})\rceil/\ell_n\big).
	}
Thus, 
	\eqn{
	(R_{\underline{e}}/\ell_n,R_{\overline{e}}/\ell_n)
	=
	\big(F^*_{n,{\sss X}}(D_{\underline{e}}+U_{e}), 
	F^*_{n,{\sss X}}(D_{\overline{e}}+U'_{e})\big) +O(1/\ell_n).
	}
By \eqref{conv-prob-degree}, the fact that 
$X_n\convd X$ and the fact that $F^*_{\sss X}$ is continuous, 
$F^*_{n,{\sss X}}(x)\convp F^*_{\sss X}(x)$ for every $x\geq 0$.
Moreover, we claim that this convergence holds
uniformly in $x$, i.e., $\sup_{x\in \R} |F^*_{n,{\sss X}}(x)-F^*_{\sss X}(x)|\convp 0.$
To see this, note that \eqref{conv-prob-degree} implies that the distribution
functions of $X_n$ and $Y_n$ converge to those of $X$ and $Y$.
Since all these random variables take on only integer values, this convergence is 
\emph{uniform}, i.e., $\sup_{k\geq 0} |F_{n,{\sss X}}(k)-F_{\sss X}(k)| \convp 0$. 
We obtain $F^*_{n,{\sss X}}$ by linearly interpolating between $F_{n,{\sss X}}(k-1)$
and $F_{n,{\sss X}}(k)$ for every $k$, so also  $F^*_{n,{\sss X}}$ converges uniformly, 
as we claimed.

By this uniform convergence, for every bounded continuous function $g\colon [0,1]^2\to \R$,
	\eqan{
	\expec_n[g(R_{\underline{e}}/\ell_n,R_{\overline{e}}/\ell_n)]
	&=
	\expec_n[g(F^*_{n,{\sss X}}(D_{\underline{e}}+U_{e}), F^*_{n,{\sss X}}(D_{\overline{e}}+U'_{e}))]\\
	&=\expec_n[g(F^*_{\sss X}(D_{\underline{e}}+U_{e}),
	F^*_{\sss X}(D_{\overline{e}}+U'_{e}))] +\op(1)\nn\\
	&=\expec_n[g(F^*_{\sss X}(X_n+U),
	F^*_{\sss X}(Y_n+U'))]+\op(1)\nn\\
	&\convp \expec[g(F^*_{\sss X}(X+U),
	F^*_{\sss X}(Y+U'))]=\expec[g(F^*_{\sss X}(X^*),
	F^*_{\sss X}(Y^*))],\nn
	}
again by \eqref{conv-prob-degree} and the fact that $(x,y)\mapsto \expec[g(F^*_{\sss X}(x+U),
F^*_{\sss X}(y+U'))]$ is continuous and bounded. Applying this to $g(x,y)=xy$, $g(x,y)=x^2$ and
$g(x,y)=y^2$ yields the required convergence. Moreover, since $F^*_{\sss X}(X^*)$ and 
$F^*_{\sss X}(Y^*)$ are uniform random variables, $\Var(F^*_{\sss X}(X^*))=\Var(F^*_{\sss X}(Y^*))=1/12.$ 
This completes the proof \ch{of convergence in (a). The equality in (a) is just \chr{\cite[Proposition~3.1]{Mesfioui2005Spearman},} see (\ref{eq:proposition3-1}).}

For part (b), we note that
	\eqn{
	\rho(G_n)=\frac{\Cov_n(X_n,Y_n)}{\Var_n(X_n)}.
	}
Since $\expec_n[X_n^2]\convp\expec[X^2]<\infty$, also
$\expec_n[X_n]\convp\expec[X]<\infty$, so that
$\Var_n(X_n)\convp \Var(X)$.
Since these limits are positive, by Slutzky's theorem,
	\eqn{
	\rho(G_n)=\frac{\Cov_n(X_n,Y_n)}{\Var(X)}(1+\op(1)).
	}
Furthermore, the random variables $(X_nY_n)_{n\geq 1}$ 
converge in distribution, and are uniformly integrable 
(since both $(X_n^2)_{n\geq 1}$ and $(Y_n^2)_{n\geq 1}$ are, which
again follows from the fact that $\expec_n[X_n^2]\convp\expec[X^2]<\infty$
and the fact that $X_n$ and $Y_n$ have the same marginals). Therefore, also
$\expec_n[X_nY_n]\convp \expec[XY]$, so that the convergence follows.
\qed

\section{Random graph examples}
\label{sec-RG-exam}
\chn{In this section} we consider four random graph models to highlight our result: the configuration model, the configuration model with intermediate vertices, the preferential attachment model and a model of complete bipartite random graphs. In Section \ref{sec-num-res-RGs}, we present the numerical results for these models.

\subsection{The configuration model}
\label{sec-CM}
The \emph{configuration model} (CM) was invented by 
Bollob\'as in \cite{Bollobas1980CM}, inspired by
\cite{Bender1978CM}. Its connectivity structure
was first studied by Molloy and Reed \cite{Molloy1995configuration,Molloy1998size}.
It was popularized by Newman, Srogatz and Watts \cite{Newman2001StrogatzWatts}, who realized
that it is a useful and simple model for real-world networks.

Given a {\it degree sequence}, namely a sequence of $n$
positive integers ${\boldsymbol d} = (d_1,d_2,\ldots, d_n)$
with $\ell_n=\sum_{i\in [n]} d_i$ assumed to be even, the configuration model (CM) on $n$
vertices and degree sequence ${\boldsymbol d}$ is constructed as follows. 
Start with $n$ vertices, \chn{labelled $1,2,\ldots,n$}, and $d_v$ half-edges adjacent to vertex $v$. The graph is constructed by randomly pairing each half-edge to some other half-edge to form an edge. Number the half-edges from $1$ to 
$\ell_n$ in some arbitrary order. Then, at each step,
two half-edges that are not already paired are chosen uniformly at random among
all the unpaired half-edges and are paired to form a single edge
in the graph. These half-edges are removed from the list of
unpaired half-edges. We continue with this procedure of choosing and pairing two 
unpaired half-edges
until all the half-edges are paired. \chn{In the resulting graph $G_n=(V_n,E_n)$ we have $|V_n|=n$, $\ell_n=2|E_n|$.} Although self-loops and double edges may occur, these become rare as $n\to\infty$ (see e.g.\ \cite{Bollobas_RG} or \cite{Janson2009CM}
for more precise results in this direction). \chn{In the analysis we keep the self-loops and 
multiple edges, so that $\ell_n=|E_n'|$. In the numerical simulation we also consider the case where the self-loops are removed, and we collapse multiple edges to
a single edge.} As we will see in the simulations, these two cases are
qualitatively similar.

We investigate the CM where the degrees are i.i.d.\ random variables, and note that
the probability that two vertices \chn{$u$ and $v$} are directly connected is close to
$d_ud_v/\ell_n$. Since this is of product form in $u$ and $v$, the degrees at either end
of an edge are close to being independent, and in fact are asymptotically independent. 
Therefore, one expects  the assortativity coefficient of
the configuration model to converge to 0 in probability, irrespective of 
the degree distribution. 

We now make this argument precise. We make the following assumptions on our degree sequence $(d_v)_{v\in V_n}$: 
\begin{Condition}[Degree regularity]
\label{cond-degree-reg}
~\\
(a) There exists a probability distribution
$(p_k)_{k\geq 0}$ such that $n_k/n\ra p_k$ for every $k\geq 1$, where $n_k=\#\{v\colon d_v=k\}$ 
denotes the number of vertices of degree $k$.\\
(b) \chn{$\expec[D_{(n)}]\ra \expec[D]$,
where $\P(D_{(n)}=k)=n_k/n$ and $\prob(D=k)=p_k$.} 
\end{Condition}

See \cite[Chapter 7]{HofstadRG}
for an extensive discussion of the CM under Condition \ref{cond-degree-reg}.

\begin{Theorem}[Convergence of \chn{the degree-degree dependency measures} for CM]
\label{thm-conv-rho-assor-CM}
Let $(G_n)_{n\geq 1}$  be a sequence of configuration models of size $n$,
for which the degree sequence $(d_v)_{v\in V_n}$ satisfies Condition \ref{cond-degree-reg}.
Then\\
	\[
	\rho^{\rm rank}(G_n)\convp 0,
	\]
and
	\[
	\rho(G_n)\convp 0.
	\]
\end{Theorem}

\proof We apply Theorem \ref{thm-conv-rho-assor}, for which we start by investigating 
\eqref{conv-prob-degree-rep}. We note that a uniform edge can be constructed by taking two
half-edges uniformly at random. Indeed, we can first draw the first half edge uniformly at random,
and this will be paired \chr{to} another half edge uniformly at random by construction of the CM.
We perform a second moment argument on $N_{k,l}=\#\{e=(u,v)\in E'_n\colon (d_u, d_v)=(k,l)\}$,
and will prove that  
	\[
	N_{k,l}/\ell_n\convp \frac{kp_k}{\expec[D]} \frac{lp_l}{\expec[D]},
	\]
For this, it suffices to prove that
	\[
	\expec[N_{k,l}]/\ell_n\ra \frac{kp_k}{\expec[D]}\frac{lp_l}{\expec[D]},
	\qquad
	\expec[N_{k,l}^2]/\ell_n^2\ra \Big(\frac{kp_k}{\expec[D]}\frac{lp_l}{\expec[D]}\Big)^2,
	\]
since then $\Var(N_{k,l}/\ell_n)=o(1).$

We note that
	\[
	\expec[N_{k,l}]=\frac{kln_kn_l}{\ell_n-1},
	\]
where $\ell_n=\sum_{v\in V_n}d_v=2|E_n|$ and $n_k=\#\{v\colon d_v=k\}$
is the number of vertices with degree $k$. Therefore, also using that $\ell_n=n\expec[D_{(n)}]$,
Condition \ref{cond-degree-reg} implies that
	\[
	\expec[N_{k,l}]/\ell_n\ra \frac{kp_k}{\expec[D]}\frac{lp_l}{\expec[D]}.
	\]
Further,
	\[
	\expec[N_{k,l}^2]/\ell_n^2=\frac{1}{\ell_n^2}\sum_{(u_1,v_1),(u_2,v_2)}
	\P(d_{u_1}=k, d_{v_1}=l, d_{u_2}=k, d_{v_2}=l).
	\]
There are four different cases, depending on $a=\#\{u_1, u_2, v_1,v_2\}$.
When $a=4$, the contribution is
	\[
	\frac{k^2 n_k(n_k-1) l^2n_l(n_l-1)}{\ell_n^2(\ell_n-1)(\ell_n-3)}=\frac{\chrr{(kn_k ln_l)^2}}{\ell_n^4}(1+O(1/n))
	\ra \Big(\frac{kp_k}{\expec[D]} \frac{lp_l}{\expec[D]}\Big)^2.
	\]
Therefore, we are left to show that the contributions due to $a\leq 3$ vanish.

When $a=3$, either one of the edges $(u_1,v_1)$ and $(u_2,v_2)$ is a self-loop,
while the other joins two other vertices (which only contributes when $k=l$), 
or both edges start in the same vertex $v$,
so that this contribution is at most
	\[
	\chr{\frac{k^2n_k(n_k-1)l^2 n_l}{\ell_n^2(\ell_n-1)(\ell_n-3)}=O(1/n)=o(1).}
	\]
When $a=2$, similar computations show that the contribution is at most $O(1/n^2)$.
When $a=1$, the edges $(u_1,v_1)$ and $(u_2,v_2)$ are self-loops
from the same vertex $v$, so that this contributes only when 
$k=l$, and then at most
	\[
	\frac{k(k-1)(k-2)(k-3)n_k}{\ell_n^2(\ell_n-1)(\ell_n-3)}=O(1/n^3)=o(1).
	\]
We conclude that \eqref{conv-prob-degree-rep} holds with 
	\[
	\P(X=k,Y=l)=\frac{kp_k}{\expec[D]} \frac{lp_l}{\expec[D]}.
	\]
In particular, $X$ and $Y$ are independent, so that $\rho^{\rm rank}=0$. This proves the first part
of Theorem~\ref{thm-conv-rho-assor-CM}.

For the second part, we note that when the degrees $(d_v)_{v\in V_n}$ are \emph{fixed},
the only random part in $\rho(G_n)$ is 
	\[
	M_n=\frac{1}{\ell_n}\sum_{e\in \chn{E'_n}} d_{\underline{e}} d_{\overline{e}}.
	\]
We perform a second moment method on this quantity.  \chr{We use that an edge $e$ 
is a pair of two specified half-edges incident to two specific vertices. Thus, we can denote $e$ by $\underline{e}=(u,s), \overline{e}=(v,t)$, where $u,v$ are the vertices to which the specific half-edges
are incident, while $s\in\{1,\ldots, d_u\}$ is the label of the half-edge 
incident to vertex $u$ and $t\in\{1,\ldots, d_v\}$ is the label of the half-edge 
incident to vertex $v$, that are paired together. The probability of pairing them together equals 
$1/(\ell_n-1)$. Therefore,}
	\[
	\expec[M_n]\chr{=\frac{1}{\ell_n}\sum_{u,v,s,t} \frac{d_{u} d_{v}}{\ell_n-1}
	=\sum_{u,v\in V_n} d_u^2d_v^2/\ell_n(\ell_n-1)= \sum_{u,v\in V_n} d_u^2d_v^2/\ell_n^2(1+O(1/n))},
	\]
where we note that we count multiple edges as frequently as they occur.
Further, \chr{and in a similar way,}
	\[	\expec[M_n^2]=(1+o(1))\sum_{u,v, u',v'\in V_n} d_u^2 d_{u'}^2 d_v^2d_{v'}^2/\ell_n^4,
	\]
so that 
\[		\frac{M_n}{\Big(\sum_{v\in V_n} d_v^2/\ell_n\Big)^2}\convp 1.
	\]
In particular, 	
	\[
	\rho(G_n)=\frac{M_n-\Big(\sum_{u,v\in V_n} d_u^2/\ell_n\Big)^2}{\sum_{u\in V_n} d_u^3/\ell_n
	-\Big(\sum_{u\in V_n} d_u^2/\ell_n\Big)^2}\convp 0,
	\]
both when $\sum_{u\in V_n} d_u^3/\ell_n\gg \Big(\sum_{u\in V_n} d_u^2/\ell_n\Big)^2$,
as well as when $\sum_{u\in V_n} d_u^3/\ell_n=\Theta \Big(\sum_{u\in V_n} d_u^2/\ell_n\Big)^2$.
\qed

\subsection{Configuration model with intermediate vertices}
\label{sec-CM-IV}

\chr{We now give an example of a strongly disassortative graph to demonstrate that $\rho(G_n)$ fails to capture obvious negative degree-degree dependencies when the degree distribution is heavy tailed. In order to do that} we adapt the configuration model slightly, by
replacing every edge by two edges that meet at a middle vertex.
Denote this graph by $\bar{G}_n=(\bar{V}_n, \bar{E}_n)$, while the 
configuration model is $G_n=(V_n,E_n)$. In this model, 
there are $n+\ell_n/2$ vertices and $\chn{|\bar{E}_n'|=}2\ell_n$ directed edges.
For $(u,v)\in \bar{E}'_n$, the degree of either vertex $u$ or vertex $v$ equals
2, and the degree of the other vertex in the edge is equal to $d_s$, where $s$ is the unique 
vertex in the original configuration model that corresponds to 
$u$ or $v$.

\begin{Theorem}[Convergence of \chn{degree-degree} dependency measures 
for CM with intermediate vertices]
\label{thm-conv-rho-assor-CMiv}
Let $(\bar{G}_n)_{n\geq 1}$  be a sequence of configuration models with intermediate vertices,
where the degree sequence $(d_v)_{v\in V_n}$ satisfies Condition \ref{cond-degree-reg}.
Then\\\ch{
	\eqn{\label{eq:spearman_cmie}
	\rho^{\rm rank}(\bar{G}_n)\convp 12\E(F_{\sss X}(X)F_{\sss X}(Y))-3= -\frac{3}{4}+3\left(\tilde{p}_1+\frac{1}{2}\tilde{p}_2\right)\left(1-\tilde{p}_1-\frac{1}{2}\tilde{p}_2\right), 
	}
where $(X,Y)=(2I+(1-I)\tilde{D}_1, 2(1-I)+I\tilde{D}_2)$ with $\tilde{D}_1, \tilde{D}_2$ i.i.d.\
random variables with $\P(\tilde{D}=k)=kp_k/\expec[D]:=\tilde{p}_k$ and $I$ an independent Bernoulli(1/2) 
random variable.} Further,
	\[
	\rho(G_n)\convp 
	\begin{cases}
	\frac{\Cov(X,Y)}{\Var(X)}&\text{if } \expec[D_{(n)}^3]\ra \expec[D^3]<\infty;\\
    	0 &\text{if } \expec[D_{(n)}^3]\ra \infty,
	\end{cases}
	\]
and, for $\expec[D_{(n)}^3]\ra \expec[D^3]<\infty$, and writing
$\mu_p=\expec[D^p]$,
	\[
	\frac{\Cov(X,Y)}{\Var(X)}=\frac{2\mu_2/\mu_1-(1+\mu_2/(2\mu_1))^2}
	{(2+\mu_3/(2\mu_1))-(1+\mu_2/(2\mu_1))^2}<0.
	\]
\end{Theorem}

The fact that the degree-degree correlation is negative is quite reasonable,
since in this model, vertices of high degree are label
only connected
to vertices of degree 2, so that there is a negative dependence
between the degrees at either end of an edge.
When $\expec[D_{(n)}^3]\ra \infty$, on the other hand,  $\rho(\bar{G}_n)\convp 0$,
which is inappropriate, as the negative dependence of 
the degrees persists.   

\proof  The first part follows directly from Theorem \ref{thm-conv-rho-assor},
since the collection of values $(\bar{d}_{\underline{e}}, \bar{d}_{\overline{e}})_{e\in \bar{E}_n'}$
only depends on the degrees $(d_v)_{v\in V_n}$ and
	\[
	\#\{e\colon \bar{d}_{\underline{e}}=l, \bar{d}_{\overline{e}}=k\}/|\bar{E}_n'|
	=(kn_k\delta_{{2},l}+ln_l\delta_{{2},k}-2n_2\indic{k=l=2})/(2\ell_n),
	\]
which converges to $\P(X=k,Y=2)$. 
Now, consider the possible values of $X$, and notice that 
\begin{align}
\label{eq:cmie_p1}
\P(X= 1)&=\tilde{p}_1/2,\\
\label{eq:cmie_p2}
\P(X=2)&=1/2+\tilde{p}_2/2,\\
\label{eq:cmie_p3}
\P(X\ge 3)&=1/2-\tilde{p}_1/2-\tilde{p}_2/2.
\end{align}
Then we obtain 
\ch{\begin{equation}
\label{eq:GX_cmie}
F^*_{\sss {X}}(\ch{x}+U)=\left\{\begin{array}{ll}\frac{1}{2}\tilde{p}_1U,&\mbox{if $\ch{x}=1$},\\
\frac{\tilde{p}_1}{2}+\Big(\frac{\tilde{p}_2}{2}+\frac{1}{2}\Big)U,&\mbox{if $\ch{x}=2$},\\
\frac{1}{2}+\sum_{k=1}^{x-1}\frac{\tilde{p}_k}{2}+\frac{\tilde{p}_x}{2}U,&\mbox{if $\ch{x}\ge 3$}.
\end{array}
\right.\end{equation}}
   Since either $X$ or $Y$ equals 2 and corresponds to the intermediate node, we further condition on $\chrr{\tilde{D}}$:
\ch{\begin{align}
\E&(F^*_{\sss X}(X^*)F^*_{\sss X}(Y^*))=\E(F^*_{\sss X}(\tilde{D}+U)F^*_{\sss X}(2+\chn{U'}))\\
&=\E(F^*_{\sss X}(2+U'))\nn\\
&\times\left[(\E(F^*_{\sss X}(1+U))\P(\tilde{D}=1)+\E(F^*_{\sss X}(2+U))\P(\tilde{D}=2)+
\E(\tilde{D}+U|\tilde{D}\ge 3)\P(\tilde{D}\ge 3)\right].\nn
\end{align}}
Now, using (\ref{eq:GX_cmie}) and substituting (\ref{eq:cmie_p1}--\ref{eq:cmie_p3}), from the last expression we readily obtain

\begin{align}
\E(F^*_{\sss X}(X^*)F^*_{\sss X}(Y^*))&=\left(\frac{\tilde{p}_1}{2}+\frac{\tilde{p}_2}{4}+\frac{1}{4}\right)\nn\\
&\qquad\times \left[\frac{1}{4}(\tilde{p}_1)^2+\Big(\frac{\tilde{p}_1}{2}+\frac{\tilde{p}_2}{4}+\frac{1}{4}\Big)\tilde{p}_2+\Big(\frac{\tilde{p}_1}{4}+\frac{\tilde{p}_2}{4}+
\frac{3}{4}\Big)(1-\tilde{p}_1-\tilde{p}_2)\right]\nn\\
&=\frac{3}{16}+\frac{1}{4}\left(\tilde{p}_1+\frac{1}{2}\tilde{p}_2\right)\left(1-\tilde{p}_1-\frac{1}{2}\tilde{p}_2\right).\nn
\end{align}
Substituting this in (\ref{eq:spearman_convergence}) and again using (\ref{eq:proposition3-1}) we obtain (\ref{eq:spearman_cmie}).

For the second part, we compute
    \chn{\[
    \frac{1}{|\bar{E}'_n|}\sum_{(u,v)\in \bar{E}_n'} \bar{d}_u\bar{d}_v
    =\frac{2}{\ell_n}\sum_{v\in V_n} d_v^2,
    \]}
and for $p\ge 2$,
    \[
    \frac{1}{|\bar{E}'_n|}\sum_{s\in \bar{V}_n} \bar{d}_s^p
    =\frac{1}{2\ell_n}2^p(\ell_n/2)+\frac{1}{2\ell_n}\sum_{v\in V_n}d_v^{p}
	=2^{p-2}+\frac{1}{2\ell_n}\sum_{v\in V_n}d_v^{p},
    \]
As a result, when $\expec[D_{(n)}^3]\ra \expec[D^3]<\infty$, we have
    \[
    \rho(\bar{G}_n)\convp \frac{2\mu_2/\mu_1-(1+\mu_2/(2\mu_1))^2}
	{(2+\mu_3/(2\mu_1))-(1+\mu_2/(2\mu_1))^2}
    <0,
      \]
      where $\mu_p=\E[D^p]$.
       \qed

\subsection{Preferential attachment model}
\label{sec-PAM} 
We discuss the general Preferential Attachment model (PAM), as formulated,
for example, in \cite[Chapter 8]{HofstadRG} or \cite[Chapter 4]{Durrett2007RG}.
The PAM is a \emph{dynamical} random graph model, and thus models
a growing network. It is defined in terms of two parameters,
$m$, which denotes the number of edges of newly added vertices, and
$\delta>-m$, which quantifies the tendency to attach to vertices
that already have a high degree. We start by defining the model for $m=1$.

We start with one \chr{vertex} having one self-loop. Suppose we 
have the graph of size $t$, which we denote by $G_t^{\sss(1)}$. \ch{Let $i$ label the vertex that appeared at time $i=1,2,\ldots$}.
Then, $G_{t+1}^{\sss(1)}$ is constructed by adding one extra vertex that 
has one edge, which forms a self-loop with probability $(1+\delta)/((2+\delta)t+1+\delta)$
and, conditionally on $G_t^{\sss(1)}$, attaches to a vertex $v\in [t]$ with probability
$(D_i(t)+\delta)/((2+\delta)t+1+\delta)$, where $D_i(t)$ is the \chn{random} degree of vertex $i$
in $G_t^{\sss(1)}$. As a result, vertices with high degree have a higher probability to be attached
to, which explains the name \emph{preferential attachment model}.

The model with $m\geq 2$ is obtained from the model with $m=1$ as follows.
Collapse vertices $m(s-1)+1,\ldots,ms$, and all of their edges,
in $(G_{t}^{\sss(1)})_{t\geq 1}$ with $\delta$ replaced by $\delta'=\delta/m$
to \chr{form} vertex $s$ in $(G_{t}^{\sss(m)})_{t\geq 1}$ with parameter $\delta$.
It is well known (see e.g., \cite{Bollobas01} where this was first derived for
$\delta=0$ and \cite[Theorem 8.3]{HofstadRG} \chr{as well as} the references in 
\cite{HofstadRG} for a more detailed literature overview) 
that the resulting graph has an asymptotic degree sequence $p_k$,
i.e.,
	\eqn{
	N_k(t)/t=\#\{i\in [t]\colon D_i(t)=k\}/t \convp p_k,
	}
where, for $k\geq m$,
    \eqn{
    \label{pkm2def}
    p_k=(2+\delta/m)\frac{\Gamma(k+\delta)\Gamma(m+2+\delta+\delta/m)}
    {\Gamma(m+\delta)\Gamma(k+3+\delta+\delta/m)}.
    }
In particular, the PAM is scale free with power-law exponent $\gamma=2+\delta/m$.
\chr{See \cite[Section 8.2]{HofstadRG} for more \chrr{details} on the scale-free behavior of the PAM.}
The next theorem investigates the behaviour of \chrr{Pearson's} correlation coefficient as well as Spearman's rho 
for the PAM:

\begin{Theorem}[Convergence of \chn{degree-degree dependency measures}
for PAM]
\label{thm-conv-rho-assor-PAM}
Let $(G_t^{\sss(m)})_{t\geq 1}$  be the PAM.
Then\\
	\eqn{
	\label{rho-rank-conv-PA}
	\rho^{\rm rank}(G_t^{\sss(m)})\convp \rho^{\rm rank}, 
	}
while
	\eqn{
	\label{rho-conv-PA}
	\rho(G_t^{\sss(m)})\convp 
	\begin{cases}
	0 &\text{if } \delta\leq m,\\
	\rho &\text{if } \delta>m,
	\end{cases}
	}
where, abbreviating $a=\delta/m$,
	\eqn{
	\label{rho-DFGM}
	\rho
	=\frac{(m-1)(a-1)[2(1+m) + a(1+3m)]}{(1+m)[2(1+m) + a(5+7m) + a^2(1+7m)].}
	}
\end{Theorem}
The value of $\rho$ \chr{in \eqref{rho-DFGM}} 
was predicted in \cite{Dorogovtsev2010correlationsPA}, and we make this analysis 
mathematically rigorous. The remainder of the section is the proof of Theorem~\ref{thm-conv-rho-assor-PAM}. It involves intermediate technical results  formulated as Lemma's~\ref{lem-PAM-Azuma}--\ref{lem-PAM-deg-deg-distrib} below.

For the PAM, it will be convenient to direct the edges from
\chr{young to old}, so that there are $mt$ \chr{directed} edges. Let $N_{k,l}(t)$ denote
the number of directed edges $e$ for which $D_{\underline{e}}(t)=k$,
$D_{\overline{e}}(t)={l}$. We will prove that there exists a probability distribution
$(q_{k,l})_{k,l\geq m}$ such that
	\eqn{
	\label{Nkl-conv}
	N_{k,l}(t)/(mt)\convp q_{k,l}.
	}
Since a uniform directed edge oriented from young to old can be obtained
by taking a uniform vertex and then a uniform edge coming out of this vertex,
this proves \eqref{conv-prob-degree-rep} with
	\eqn{
	\label{pkl-def-PAM}
	p_{kl}=\P(X=k,Y=l)=\tfrac{1}{2}(q_{k,l}+q_{l,k}).
	}
\chr{In particular, by Theorem \ref{thm-conv-rho-assor}(a), this proves \eqref{rho-rank-conv-PA}
in Theorem \ref{thm-conv-rho-assor-PAM}.}
We follow the proof of \cite[Theorem 8.2]{HofstadRG}, which, in turn,
is strongly inspired by the proof in \cite{Bollobas01}.

Proofs for convergence of the degree sequence typically 
consist of two key steps. The first is a martingale concentration 
argument in Lemma~\ref{lem-PAM-Azuma}. 

\begin{Lemma}[Convergence of degree-degree counts]
\label{lem-PAM-Azuma}
For every $k,l$, there exists a $C>0$ such that,
	\eqn{
	\label{aim-AH}
	\prob\Big(\max_{k,l} |N_{kl}(t)-\expec[N_{kl}(t)]|\geq C\sqrt{t\log{t}}\Big)=o(1).
    	}
    	\end{Lemma}

\proof The proof for the degree distribution in \cite{HofstadRG} 
applies almost verbatim (see, in particular, \cite[Proposition 8.4]{HofstadRG}
and its proof). Indeed, the proof relies on a martingale argument. Define the Doob-martingale,
\chr{for $t=0, \ldots, n$},
	\[
	M_n=\expec[N_{kl}(t)\mid G_{n}^{\sss(m)}].
	\]
The crucial observation is that $(M_n)_{n=0}^{t}$ is  a martingale with $M_t=N_{kl}(t)$ and
$M_0=\expec[N_{kl}(t)]$ \chr{that satisfies}
	\eqn{
	\label{hoeffding}
	|M_n-M_{n-1}|\leq 4m.
	}
\chr{We prove \eqref{hoeffding} below.} The Azuma-Hoeffding inequality 
\cite{Azuma1967, Hoeffding1963} then proves 
\eqref{aim-AH} for any $C>4[4m]^2$. \chr{Indeed, 
	\eqan{\nn
	\prob\Big(|N_{kl}(t)-\expec[N_{kl}(t)]|\geq A\Big)
	&=\prob\Big(|M_t-M_0|\geq A\Big)\leq \e^{-A^2/(2t[4m]^2)}.
	}
Taking $A=C\sqrt{t\log{t}}$ with $C^2>4[4m]^2$ proves that 
	\[
	\prob\Big(|N_{kl}(t)-\expec[N_{kl}(t)]|\geq C\sqrt{t\log{t}}\Big)=o(1/t^2),
	\]
so that even 
	\eqan{\nn
	&\prob\Big(\max_{k,l} |N_{kl}(t)-\expec[N_{kl}(t)]|\geq C\sqrt{t\log{t}}\Big)\\
	&\qquad\leq (m t)^2 \max_{k,l}\prob\Big(\max_{k,l} |N_{kl}(t)-\expec[N_{kl}(t)]|\geq C\sqrt{t\log{t}}\Big)=o(1).\nn
	}
This completes the proof of Lemma \ref{lem-PAM-Azuma} assuming \eqref{hoeffding}.
}

\chr{We complete the proof by deriving \eqref{hoeffding}. For this, it will be convenient to introduce some 
further notation. Let $e\in[mt]$ label the edges. Let $v_e=\lceil e/m\rceil$ denote the vertex from which 
the $e$th edge emanates, and $V_e$ (which is a random variable) the vertex to which the $e$th edge points. Then,
	\[
	N_{k,l}(t)=\sum_{e\in[mt]} \indic{D_{v_e}(t)=k, D_{V_e}(t)=l}.
	\]
As a result,
	\[
	M_n-M_{n-1}=\sum_{e\in[mt]} \Big[\prob(D_{v_e}(t)=k, D_{V_e}(t)=l\mid G_{n})-
	\prob(D_{v_e}(t)=k, D_{V_e}(t)=l\mid G_{n-1})\Big],
	\]
where we abbreviate $G_{n}=G_{n}^{\sss(m)}$. We let $(G_l')_{l\geq 0}$ denote the PAM with 
$G_{n-1}'=G_{n-1}$, while the evolution of $(G_l')_{l\geq 0}$ after time $n-1$ is the same in distribution as that of $(G_l)_{l\geq 0}$, but conditionally independent of it given $G_{n-1}=G_{n-1}'$. 
Let $D_i'(t)$ denote the degree of vertex $i$ in $G_t'$. Then, 
	\eqan{\nn
	\prob(D_{v_e}(t)=k, D_{V_e}(t)=l\mid G_{n-1})
	&=\prob(D_{v_e}'(t)=k, D_{V_e}'(t)=l\mid G_{n-1})\\
	&=\prob(D_{v_e}'(t)=k, D_{V_e}'(t)=l\mid G_{n-1}, G_n\setminus G_{n-1}),\nn
	}
where $G_n\setminus G_{n-1}$ is shorthand for the edges of $G_n$ that are not in $G_{n-1}$.
The last step is due to the conditional independence of the evolution after time $n-1$ in $(G_t')_{t\geq 0}$. 
Thus,
	\[
	\prob(D_{v_e}(t)=k, D_{V_e}(t)=l\mid G_{n-1})
	=\prob(D_{v_e}'(t)=k, D_{V_e}'(t)=l\mid G_{n}).
	\]
We conclude that
	\eqan{\nn
	M_n-M_{n-1}&=\sum_{e\in[mt]} \Big[\prob(D_{v_e}(t)=k, D_{V_e}(t)=l\mid G_{n})-
	\prob(D_{v_e}'(t)=k, D_{V_e}'(t)=l\mid G_{n})\Big].
	}	
When $V_{e}>n$, clearly $\prob(D_{v_e}(t)=k, D_{V_e}(t)=l\mid G_{n})=
\prob(D_{v_e}'(t)=k, D_{V_e}'(t)=l\mid G_{n})$, as the degrees of vertices $i$ with $i>n$ are independent of $G_n$. Thus, we can restrict to $V_e\leq n$. Further, when $v_e>n$, then 
$D_{v_e}(t)$ is independent of $G_{n}$, so that
	\eqan{\nn
	&\prob(D_{v_e}(t)=k, D_{V_e}(t)=l\mid G_{n})-
	\prob(D_{v_e}'(t)=k, D_{V_e}'(t)=l\mid G_{n})\\
	&\quad=\prob(D_{v_e}(t)=k)\Big[\prob(D_{V_e}(t)=l\mid G_{n})-\prob(D_{V_e}'(t)=l\mid G_{n})\Big].\nn
	}
Note that $D_{V_e}(n-1)=D_{V_e}'(n-1)$ a.s., $\prob(D_{V_e}(t)=l\mid G_{n}, D_{V_e}(n)=j)
=\prob(D_{V_e}(t)=l\mid D_{V_e}(n)=j)$, and
	\[
	\prob(D_{V_e}'(t)=l\mid G_{n},D_{V_e}'(n)=j)
	=\prob(D_{V_e}'(t)=l\mid D_{V_e}'(n)=j)=\prob(D_{V_e}(t)=l\mid D_{V_e}(n)=j).
	\]
Thus, using that
	\eqan{\nn
	\prob(D_{V_e}(t)=l\mid G_{n})&=\expec[\prob(D_{V_e}'(t)=l\mid D_{V_e}(n))\mid G_n],\\
	\nn\prob(D_{V_e}'(t)=l\mid G_{n})&=\expec[\prob(D_{V_e}'(t)=l\mid D_{V_e}'(n))\mid G_n],
	}
we obtain at
	\[
	|\prob(D_{V_e}'(t)=l\mid D_{V_e}(n))-\prob(D_{V_e}'(t)=l\mid D_{V_e}'(n))|
	\leq \indic{D_{V_e}(n)\neq D_{V_e}'(n)}.
	\]
Taking expectations \chrr{yields}
	\[
	\Big|\prob(D_{v_e}(t)=k, D_{V_e}(t)=l\mid G_{n})-
	\prob(D_{v_e}'(t)=k, D_{V_e}'(t)=l\mid G_{n})\Big|
	\leq \prob(D_{V_e}(n)\neq D_{V_e}'(n)\mid G_{n}).
	\]
In a similar way, we see that for $v_e\leq n$,
	\eqan{\nn
	&|\prob(D_{v_e}(t)=k, D_{V_e}(t)=l\mid G_{n})-
	\prob(D_{v_e}'(t)=k, D_{V_e}'(t)=l\mid G_{n})|\\
	&\quad\leq \prob(D_{V_e}(n)\neq D_{V_e}'(n)\mid G_{n})
	+\prob(D_{v_e}(n)\neq D_{v_e}'(n)\mid G_{n}).\nn
	}
We conclude that
	\eqan{\nn
	|M_n-M_{n-1}|
	&\leq \sum_{e\in[mt]} \big[\prob(D_{V_e}(n)\neq D_{V_e}'(n)\mid G_{n})
	+\prob(D_{v_e}(n)\neq D_{v_e}'(n)\mid G_{n})\big]\leq 4m.
	}
}
\vskip-1.2cm\hfill \qed

\vskip0.6cm
\chr{We continue with the proof of \eqref{Nkl-conv}.}
The second key step \chr{the proof of \eqref{Nkl-conv}} is to prove that, for each $k,l$,
	\eqn{
	\label{mean-Nkl-conv}
	\lim_{t\ra \infty} \expec[N_{kl}(t)]/(mt) =q_{k,l}.
	}
We \chr{sum over} the vertex $s$ that has degree $l$ at time $t$, \chr{and condition on} the 
degree $r\geq m$ of the vertex to which the edge of vertex $s$ is attached.
This yields
	\eqn{
	\label{Nkl-mean-a}
	\expec[N_{kl}(t)]=
	m\sum_{s=1}^t \sum_{r\geq m} 
	\frac{(r+\delta)}{(2m+\delta)s}\expec[N_r(s)] 
	\Big[\prob\big(\chr{B_{r+1}}[s+1,t]=k, \chr{B_{m}}[s+1,t]=l\big)
	+O(1/s)\Big],
	}
where $\chr{B_{m}}[s+1,t]$ is $m$ plus the number of edges attached to vertex $s$ between time $s+1$ and $t$,
while $\chr{B_{r+1}}[s+1,t]$ is $r$ plus the number of further edges attached to the vertex of
degree $r$ to which the edge of vertex $s$ is attached. The $\chr{O(1/s)}$ term is due to 
contributions where at least \emph{two} edges of vertex $s$ are attached to the same vertex
of degree $r$, \chr{and also due to the fact that the probability of attaching the $j$th edge of vertex 
$s$ to a vertex of degree $r$ at time $s$ is actually equal to $\frac{(r+\delta)}{(2m+\delta)s+(j-1)(2+\delta/m)+1+\delta/m}$, which is $\frac{(r+\delta)}{(2m+\delta)s}(1+O(1/s)).$} Further,
	\[
	\prob\big(\chr{B_{r+1}}[s+1,t]=k, \chr{B_{m}}[s+1,t]=l\big)
	=\prob(\chr{B_{r+1}}[s+1,t]=k)\prob(\chr{B_{m}}[s+1,t]=l)+O(1/t),
	\]
since the dependence between the two probabilities is entirely due to the 
fact that edges that contribute to $\chr{B_{r+1}}[s+1,t]$ cannot contribute to 
$\chr{B_{m}}[s+1,t]$. Indeed, \chr{$(B_{r+1}[s+1,t], B_{m}[s+1,t])$ is equal in distribution 
to the number of balls in two urns at time $m(t-s)$, where we start with $r+1$ and $m$ 
balls at time $0$, and in each draw, we draw a ball in each of the urns with probability 
equal to the number of balls plus $\delta$ and then replace it with two balls. 
Knowing how many balls are put into the first urn only gives us 
information about how many balls cannot be put into the second urn, 
so the balls in the different urns are close to independent.}
We study these probabilities \chrr{now}:

\begin{Lemma}[Growth of degrees in PAM]
\label{lem-PAM-degrees}
For all $k\geq r\geq m$ and $a\in (0,1)$, 
	\[
	\lim_{s\ra \infty} \prob(\chr{B_{r}}[as,s]=k)=P_k(a;r),
	\]
where, for each $r\geq m$ and $a\in (0,1)$,  $(P_k(a;r))_{k\geq r}$ is a probability measure.
\end{Lemma}

\proof We note that $(\chr{B_{r}}[s,ts])_{t\geq 1}\convd (Z_t)_{t\geq 1}$, as $s\rightarrow \infty$,
where $(Z_t)_{t\geq 0}$ is a pure birth process, which increases by 1 at rate 
$m(Z_t+\delta)/((2m+\delta)t)$ at time $t$. Indeed, when $\chr{B_{r}}[s,ts]=k$,
then each of the $m$ edges of vertex $st+1$ has probability $(k+\delta)/[(2m+\delta)(st)] +O(1/s^2)$
of being attached to the vertex that has degree $k$ at time $ts$, and thus of increasing
$\chr{B_{r}}[s,ts]$ to $k+1$. Thus, within a short time interval $[t, t+dt]$ and conditionally on
$\chr{B_{r}}[s,ts]=k$, the probability that $\chr{B_{r}}[s,(t+dt)s]=k+1$ is equal to 	
	\[	\chr{s dt \Big[m(k+\delta)/[(2m+\delta)(st)] +O(1/s^2)+o(1)\Big]
	\rightarrow dt \frac{m(k+\delta)}{(2m+\delta)t}+o(\chrr{dt}),}
	\]
\chr{as $s\rightarrow \infty$. This is the birth rate of the pure birth process $(Z_t)_{t\geq 1}$.}

We next study the limiting birth process, for which is it useful to make a time change.
With $b_t=Z_{\e^{(2+\delta/m)t}}$, $(b_t)_{t \geq 0}$ is a birth process 
that grows at rate $b_t+\delta$ at time $t$. Define
	\[
	f_{r,k}(t)=\prob(b_t=k\mid b_0=r).
	\]
Then,
	\[
	\frac{\partial}{\partial t}  f_{r,k}(t)=-(k+\delta) f_{r,k}(t)+\chr{(k-1+\delta)}f_{r,k-1}(t).
	\]
This set of differential equations is solved by $f_{r,r}(t)=\e^{-(r+\delta)t}$ and,
for $k\geq r+1$,
	\[
	f_{r,k}(t)=(k-1+\delta) \e^{-(k+\delta)t} \int_0^t \e^{(k+\delta)s} f_{r,k-1}(s)ds.
	\]
This can be solved by, for $k\geq r+1$,
	\[
	f_{r,k}(t)=\prob(b_t=i\mid b_0=r)=\frac{\Gamma(k+\delta)}{\Gamma(r+\delta)} \e^{-(k+\delta)t}
	\sum_{j=0}^{k-r} \alpha_{j,k}\e^{jt},
	\]
where $\alpha_{0,k}=-\sum_{j=0}^{k-1} \alpha_{j,k-1}/(j+1)$, while, for $j\geq 1$,
	\[
	\alpha_{j,k}=\alpha_{j-1,k-1}/j.
	\]
As a result, for all $a\in (0,1)$,
	\[
	\lim_{t\ra \infty} \prob(\chr{B_{r}}[at,t]=k)=\prob(Z_{1/a}=k\mid Z_1=r)=f_{r,k}((2+\delta/m)^{-1}\log(1/a)).
	\]
\chr{Note that $P_r(a;r)$ is the probability that the birth process has no births.} We thus compute that $P_r(a;r)=f_{r,r}((2+\delta/m)^{-1}\log(1/a))=a^{(r+\delta)/(2+\delta/m)}$ for $k=r$, 
while
	\[
	P_k(a;r)=f_{r,k}((2+\delta/m)^{-1}\log(1/a))
	=\frac{\Gamma(k+\delta)}{\Gamma(r+\delta)} a^{(k+\delta)/(2+\delta/m)}
	\sum_{j=0}^{k-r} \alpha_{j,k}a^{-j/(2+\delta/m)}.
	\]
\vskip-1.18cm \hfill\qed
\vskip0.6cm

\vskip0.5cm

\noindent
We continue from \eqref{Nkl-mean-a}, and rewrite it as
	\eqn{
	\label{Nkl-mean-b}
	\expec[N_{kl}(t)]/(mt)=
	\chr{\sum_{r\geq m} \expec\Big[\frac{(r+\delta)}{(2m+\delta)Ut}\expec[N_r(Ut)] \prob(\chr{B_{r+1}}[Ut,t]=k\mid U)\prob(\chr{B_{m}}[Ut,t]=l\mid U)\Big]
	+O(\chr{\log{t}}/t),}
	}
where $U$ has a uniform distribution, we interpret $Ut=\lceil Ut\rceil$, \chr{and the outer expectation is 
over $U$ only.} Using that $\expec[N_r(s)]/s=p_r+O(1/s)$ (see \cite[Proposition 8.4]{HofstadRG}), we further arrive at
	\eqn{
	\label{Nkl-mean-c}
	\expec[N_{kl}(t)]/(mt)=
	\sum_{r\geq m} 
	\frac{r+\delta}{2m+\delta} p_r \chr{\expec\Big[\prob(\chr{B_{r+1}}[Ut,t]=k\mid U)\prob(\chr{B_{m}}[Ut,t]=l\mid U)\Big]}
	+o(1).
	}
By Lemma \ref{lem-PAM-degrees}, this converges to
	\eqn{
	\label{Nkl-mean-d}
	\expec[N_{kl}(t)]/(mt)\ra 
	q_{k,l}\equiv \sum_{r\geq m} 
	\frac{r+\delta}{2m+\delta} p_r \expec[P_k(U;r)P_l(U;m)].
	}
This proves \eqref{mean-Nkl-conv}, and \chr{thus, by Theorem \ref{thm-conv-rho-assor}(a),}
proves the convergence of the rank correlation in \eqref{rho-rank-conv-PA} \chr{in Theorem \ref{thm-conv-rho-assor-PAM}}.
\medskip

For the convergence of the correlation coefficient in 
\eqref{rho-conv-PA} \chr{in Theorem \ref{thm-conv-rho-assor-PAM}}, we \chr{aim to use 
Theorem \ref{thm-conv-rho-assor}(b) and thus start by investigating the convergence of moments of $X_n$.
By \eqref{relation-moment-degrees}, and letting $\expec_n$ denote the conditional expectation given $G_n$,
	\[
	\expec_n[X_n^2]=\frac{1}{n} \sum_{i\in[n]} D_i(n)^3.
	\]
Thus, we are lead to studying sums of powers of degrees.}
To analyze the limit of sums of powers of degrees, we rely on the following lemma:

\begin{Lemma}[Sum of powers of degrees in PAM]
\label{lem-PAM-sum-degrees}
For all \chr{$p<\chrr{\gamma=2+\delta/m},$}
	\[
	\frac{1}{n} \sum_{\chn{i\in [n]}} D_i(n)^p \convp \mu_p=\sum_{k\geq m} k^p p_k<\infty.
	\]
\end{Lemma}

\proof We note that $\sum_{i\in [n]} D_i(n)^p=\sum_{k\geq m} k^p N_k(n)$. Under the conditions
stated,  for every $k_n\rightarrow \infty$,
	\[
	\sum_{k\geq m} k^p N_k(n)=\sum_{m\leq k\leq k_n} k^p N_k(n)+\op(n).
	\]
This follows since, for any $\vep>0$, \chr{$k> k_n$ implies that $k^{\vep}/k_n^{\vep}>1$, so that}
	\[
	\sum_{k>k_n} k^p N_k(n)\leq k_n^{-\vep}\sum_{k\geq m} k^{p+\vep} N_k(n)
	=k_n^{-\vep} \frac{1}{n} \sum_{i\in [n]} D_i(n)^{p+\vep}.
	\]
By the analysis in \cite[Section 8.1 and 8.6]{HofstadRG}, when \chr{$p+\vep<\gamma+1=3+\delta/m,$}
	\[
	\limsup_{n\rightarrow \infty} \expec\Big[\frac{1}{\chn{n}}\sum_{i\in [n]} D_i(n)^{p+\vep}\Big]<\infty.
	\]
Therefore, by the Markov inequality, $\sum_{k>k_n} k^p N_k(n)=\op(n)$.

Now, since $\max_{k} |N_k(n)-p_k|\leq \sqrt{Cn\log{n}}$ whp by \cite[Proposition 8.4]{HofstadRG},
	\[
	\sum_{m\leq k\leq k_n} k^p N_k(n)=t \sum_{m\leq k\leq k_n} k^p p_k + \Op(k_n^{p+1}\sqrt{n\log{n}}).
	\]
This proves the claim.
\qed
\vskip0.5cm

\noindent
\chr{It follows from Lemma \ref{lem-PAM-sum-degrees} \chr{that for $3<\gamma=2+\delta/m,$}
	\[
	\expec_n[X_n^2]=\frac{1}{n}\sum_{i\in \chn{[n]}} D_i(n)^3=
	B(1+\op(1)).
	\]
where $B$ is a constant. As a result,
	\eqn{
	\label{rho-PAM-1}
	\rho(G_\chn{n})\convas \rho=\Cov(X,Y)/\Var(X)
	=\frac{\sum_{k,l} kl q_{k,l}-\expec[X]^2}
	{\expec[X^2]-\expec[X]^2}.
	}
This proves \eqref{rho-conv-PA} in Theorem \ref{thm-conv-rho-assor-PAM}
when $\delta>m$. For $\gamma<3$, instead, $D_1(\chn{n})/n^{1/\gamma}\convas \xi$, for some strictly positive 
random variable $\xi$ (see e.g., \cite[Sections 8.1 and 8.6]{HofstadRG}). Therefore, $\expec_n[X_n^2]\geq \xi^3 n^{3/\gamma-1}(1+o(1))$. Further, the majority of edges of high degree vertices is young, so that 
	\eqn{
	\label{cross-product}
	\expec_n[X_nY_n]=\op(n^{3/\gamma-1}).
	}
Indeed, fix $T_n$ such that $T_n\rightarrow \infty$ and $T_n=o(n)$. There are at most $mT_n$ edges between vertices with index at most $T_n$, and, since the maximal degree is $\Op(n^{1/\gamma})$, these contribute at most $\Op(n^{2/\gamma-1} T_n)$. For the other edges, one of the vertices involved 
was born after time $T_n$. Since $\max_{i\geq T_n}D_i(n)=\op(n^{1/\gamma})$,
the contribution of these edges is at most
	\[
	\op(n^{1/\gamma})\expec_n[X_n+Y_n].
	\]
In turn, $\expec_n[X_n+Y_n]=\Op(n^{(2/\gamma-1)\wedge 1})$, which completes the proof of \eqref{cross-product}. This implies that $\rho(G_n)\convp 0$, which proves \eqref{rho-conv-PA} in Theorem \ref{thm-conv-rho-assor-PAM} when $\delta<m$. For $\delta=m$, so that $\gamma=3$, $\sum_{i\in \chn{[n]}} D_i(n)^3=\Theta_{\sss \prob}(n\log{n})(1+\op(1))$. As a result, also in this case
$\rho(G_n)\convas 0$ for $\delta\leq m$.
}
\qed
\bigskip

\noindent
\chr{We continue with the proof of \eqref{rho-DFGM} in Theorem \ref{thm-conv-rho-assor-PAM}.}
To compute expectations involving $X$, we often rely on the following lemma:

\begin{Lemma}[Degree on one side of uniform edge]
\label{lem-degree-one-side-edge}
For every function $f\colon \N\to \R$,
	\[
	\expec[f(X)]=\sum_{k\geq m} f(k)\frac{kp_k}{2m}.
	\]
\end{Lemma}

\proof Let $f$ be bounded, and let $X_n$ be the degree at the bottom of a 
uniform edge. Then, 
	\eqan{\nn
	\expec[f(X_n)\mid G_n^{\sss(m)}]
	&=\frac{1}{|E_n'|} \sum_{e\in E_n'} f(D_{\underline{e}}(n))
	=\frac{1}{2mn}\sum_{v\in [n]} f(D_v(n)) D_v(n)
	=\frac{1}{2m}\sum_{k\geq m} f(k)k N_k(n)/n.
	}
Taking the limit of $n\ra \infty$ and using that $N_k(n)/n\convp p_k$,
as well as $X_n\convd X$ proves the claim. 
\qed

\medskip
Lemma \ref{lem-degree-one-side-edge} allows us to identify the r.h.s.\ of \eqref{rho-PAM-1} as 
	\[
	\rho=\Cov(X,Y)/\Var(X)
	=\frac{(2m)^2\sum_{k,l} klq_{k,l} -\lambda_2^2}{2m\lambda_3-\lambda_2^2},
	\]
where $\lambda_a=\sum_{k\geq m} k^a p_k$. To identify the limit, we 
follow \cite{Dorogovtsev2010correlationsPA}. \chr{Recall the definition of $p_{kl}$ in \eqref{pkl-def-PAM}.}

\begin{Lemma}[Asymptotic degree-degree distribution for PAM]
\label{lem-PAM-deg-deg-distrib}
For all $k,l\geq m$,
	\eqan{
	\label{akl-formula}
	&\chr{p_{kl}=\P(X=k,Y=l)}\\
	&\quad=(2+\delta/m)\frac{\Gamma(m+2+\delta+\delta/m)}{\Gamma(m+\delta)^2}
	\frac{\Gamma(l+\delta)\Gamma(k+\delta)}{\Gamma(k+2+\delta)\Gamma(l+k+2+2\delta+\delta/m)}\nn\\
	&\qquad \times
	\chrr{\Big[ \sum_{j=m+1}^k{{k+l-j-m}\choose {l-m}} {{j+k+2+2\delta+\delta/m}\choose {k+1+\delta}}+\sum_{j=m+1}^l {{k+l-j-m}\choose {k-m}} {{j+l+2+2\delta+\delta/m}\choose {l+1+\delta}}\Big]}.\nn
	}
Consequently, \eqref{rho-DFGM} follows.
\end{Lemma}
\proof To compute $\prob(X=k,Y= l)$, we let $M_{kl}(t)$ denote the number of edges at time $t$ where one side has degree $k$ and the other side degree $l$,
so that
	\[
	p_{kl}=\lim_{t\ra \infty} \expec[M_{kl}(t)]/(2mt).
	\]
We note that $M_{kl}(t)$ satisfies the recursion relation
	\eqan{\nn
	\expec[M_{kl}(t+1)]-\expec[M_{kl}(t)]&=m\frac{(k\vee l)-1+\delta}{(2m+\delta)t} \expec[N_{k\vee l-1}(t)]
	\indic{k\wedge l=m}\\
	&\quad  +m\frac{k-1+\delta}{(2m+\delta)t} \expec[M_{k-1,l}(t)]+m\frac{l-1+\delta}{(2m+\delta)t} \expec[M_{k,l-1}(t)]\nn\\
	&\quad -m\frac{k+\delta}{(2m+\delta)t} \expec[M_{k,l}(t)]-
	m\frac{l+\delta}{(2m+\delta)t} \expec[M_{k,l}(t)] +O(1/t^2).\nn
	}
\chr{It is not clear that the left-hand side converges since we only know that $\expec[M_{k,l}(t)]/(2mt)\rightarrow p_{kl}$, and we will show this now.} Indeed, since $\expec[M_{k,l}(t)]/(2mt)\ra p_{kl}$ and 
$\expec[N_{k}(t)]/t\ra p_{k}$, we arrive at the claim that,
for all $k,l$ with $k\vee l\geq m+1$,
	\eqan{\nn
	&\lim_{t\ra \infty} \expec[M_{kl}(t+1)]-\expec[M_{kl}(t)]\\
	&\quad=2m^2\frac{(k\vee l)-1+\delta}{2m+\delta}p_{k-1}
	\indic{k\wedge l=m}+
	2m^2\frac{k-1+\delta}{2m+\delta} p_{k-1,l}+2m^2\frac{l-1+\delta}{2m+\delta} p_{k,l-1}
	-2m^2\frac{k+l+2\delta}{2m+\delta} p_{k,l}.\nn
	}
Since $\lim_{t\ra \infty} \expec[M_{kl}(t)]/(2mt)=p_{kl}$,
we must therefore have that $\lim_{t\ra \infty} \expec[M_{kl}(t+1)]-\expec[M_{kl}(t)]=2mp_{kl}$, so that
	\[
	p_{kl}=m\frac{(k\vee l)-1+\delta}{2m+\delta}p_{k\vee l-1}
	\indic{k\wedge l=m}+
	m\frac{k-1+\delta}{2m+\delta} p_{k-1,l}+m\frac{l-1+\delta}{2m+\delta} p_{k,l-1}
	-m\frac{k+l+2\delta}{2m+\delta} p_{k,l},
	\]
and
	\eqn{
	\label{pkl-recur-DorFerGolMen}
	(k+l+2+2\delta+\delta/m)p_{kl}=((k\vee l)-1+\delta)p_{k\vee l-1}
	\indic{k\wedge l=m}+
	(k-1+\delta)p_{k-1,l}+(l-1+\delta)p_{k,l-1}.
	}
This is equivalent to \cite[(12)]{Dorogovtsev2010correlationsPA}.
This can be worked out to yield
	\eqan{\nn
	p_{kl}&=\sum_{j=m+1}^k {{\chrr{k+l-j-m}}\choose {k-j}} \frac{\Gamma(k+\delta)}{\Gamma(\chrr{j-1}+\delta)}
	\frac{\Gamma(\chrr{l}+\delta)}{\Gamma(m+\delta)} \frac{\Gamma(j+k+2+2\delta+\delta/m)}{\Gamma(l+k+3+2\delta+\delta/m)}p_{j-1}\\
	&\quad +\sum_{j=m+1}^l {{k+l-j-m}\choose {l-j}} \frac{\Gamma(k+\delta)}{\Gamma(\chrr{j-1}+\delta)}
	\frac{\Gamma(\chrr{l}+\delta)}{\Gamma(\chrr{m}+\delta)} \frac{\Gamma(j+\chrr{l}+2+2\delta+\delta/m)}{\Gamma(l+k+3+2\delta+\delta/m)}p_{j-1}.\nn
	}
Substituting \eqref{pkm2def}, we arrive at \eqref{akl-formula}. 

The computation to go from \eqref{pkl-recur-DorFerGolMen} to \eqref{rho-DFGM}
is performed in \cite[(12)]{Dorogovtsev2010correlationsPA}, and applies verbatim.
%
%
\qed

\subsection{Asymptotically random \chn{Pearson's coefficient}: collection of complete bipartite graphs}
\label{sec:random_assortativity}

In this section, we present an example where $\rho(G_n)$ in (\ref{eq:rhoG}) converges to a random 
variable when the number of vertices tends to infinity. 
\chn{For $|V_n|=n$,} under the assumptions of Theorem~\ref{thm-corr-coef-pos-scale-free}, we have
	\chn{\begin{align}
	\label{eq:crossproducts1}
	\sum_{(u,v)\in E_n'} D_uD_v&
	\le \max_{v\in V_n}d_v\sum_{(u,v)\in E'_n} D_u
	=\max_{v\in V_n}D_v\Big(\sum_{v\in V_n} D^2_v\Big)
	\le C^2n^{1/\gamma+(2/\gamma\vee 1)},\\
	\label{eq:crossproducts2}
	\sum_{(u,v)\in E_n'} D_uD_v&\ge \max_{v\in V_n}D_v\ge c n^{1/\gamma},\\
	\label{eq:crossproducts3}
	\sum_{(u,v)\in E'_n} D_uD_v
	&\ge\sum_{v\in V_n} D^2_v\ge  cn^{2/\gamma\vee 1}.
	\end{align}}
Further, from the proof of Theorem~\ref{thm-corr-coef-pos-scale-free}, we know that
	\eqan{
	\label{eq:cubes}\sum_{v\in V_n} D_v^3&\geq (\max_{v\in V_n} D_v)^3
	\geq c^3 n^{3/\gamma},
	}
and
	\eqan{
	\label{eq:squares}
	\frac{1}{|E_n'|}\Big(\sum_{v\in V_n} D_v^2\Big)^2
	&\leq (C^2/c) n^{(4/\gamma-1)\vee 1},
	}
where we see that (\ref{eq:squares}) is vanishing compared to (\ref{eq:cubes}). 
The convergence of (\ref{eq:rhoG}) to a random variable can only take place if 
the crossproducts on the left-hand side of (\ref{eq:crossproducts1} -- \ref{eq:crossproducts3}) 
are of the same order of magnitude as the left-hand side of (\ref{eq:cubes}). 
As we see from the above, this is possible for $\gamma\in (1,3)$. 

Below we present an example where $\rho(G_n)$ indeed 
converges to a random variable. However, 
due to slow convergence, a substantially larger computational 
capacity is needed in order to \chr{approximate} the limiting distribution.

Take $((X_i,Y_i))_{i=1}^n$ to be an i.i.d.\ sample of integer random
variables as in \eqref{eq:X}, where $\alpha_1=\alpha_2=\beta_1=b$,
$\beta_2=ab$ for some $b>0$ and $a>1$.
Then, for $i=1, \ldots, n$, we create a complete bipartite graph of 
$X_i$ and $Y_i$ vertices, respectively. These $n$ complete bipartite graphs
are not connected to one another. We denote such a collection of $n$ bipartite 
graphs by $G_n$. The graph $G_n$  has
\chn{$|V_n|=\sum_{i=1}^n (X_i+Y_i)$ vertices and $|E'_n|=2\sum_{i=1}^n X_iY_i$ directed edges.}  
Further, \chn{if $D_v$ denotes the random degree of vertex $v$, then we obtain}
    	\[
    	\sum_{v\in V_n} D_v^p
    	=\sum_{i=1}^n (X_i^pY_i+Y_i^pX_i),
	\qquad
	\sum_{(u,v)\in E'_n}D_uD_v=2\sum_{i=1}^n (X_iY_i)^2.
    	\]
Assume that the $\chrr{\xi_j}$'s in \eqref{eq:X} satisfy \eqref{tail-U} with $\gamma\in(2,4)$, so that $\expec[\chrr{\xi}^2]<\infty$,
but $\expec[\chrr{\xi}^4]=\infty$. As a result, $|E'_n|/n \convp 2\expec[XY]<\infty$
and $\frac{1}{n} \sum_{v\in V} D_v^2\convp \expec[XY(X+Y)]<\infty$ when $\gamma\in (3,4)$, while,
for $\gamma\in(2,3)$,
	\eqn{
	n^{-3/\gamma}\sum_{v\in V} D_v^2
	=n^{-3/\gamma}\sum_{i=1}^n (X_i^2Y_i+Y_i^2X_i)
	\convd Z,
	}
for some random variable $Z$. [For $\gamma=3$, this sum grows as a slowly varying function in $n$, but this case is very similar and will thus be omitted.]
Further,
	\[
	n^{-4/\gamma}b^{-4}\sum_{i=1}^n (X_i^3Y_i+Y_i^3X_i)
	\convd (a^3+a) Z_1+ 2Z_2,
	\qquad
	n^{-4/\gamma}b^{-4}\sum_{i=1}^N (X_iY_i)^2
	\convd a^2 Z_1+Z_2,
	\]
where $Z_1$ and $Z_2$ and two independent stable distributions with parameter $\gamma/4$.
\chr{Therefore, using \refeq{eq:rhoG} and the fact that $4/\gamma >(6/\gamma-1)\wedge 1$, we arrive at}
	\[
	\rho(G_n)\convd \frac{2a^2Z_1+2Z_2}{(a+a^3)Z_1+2Z_2},\; \mbox{as $n\to\infty$}.
	\]
which is a proper random variable taking values in $(2a/(1+a^2), 1)$.

For convergence of the rank correlation, we note that 
	\[
	\prob(X_n=k,Y_n=l)\rightarrow \prob(X=k,Y=l)=\frac{kl}{\expec[X_1Y_1]} \prob(X_1=k, Y_1=l),
	\]
where we recall that $(X_1,Y_1)$ is as in \eqref{eq:X}, while $(X,Y)$ are the degrees at either side
of a uniformly chosen edge. Thus, convergence of the rank correlation follows from Theorem \ref{thm-conv-rho-assor}(a).

\section{Numerical results} 
\label{sec-num-res-RGs}
In this section, we present numerical examples that illustrate our results. 

\subsection{Numerical results for configuration models and preferential attachment model} We have generated random graphs of different sizes using the configuration model in Section~\ref{sec-CM}, the configuration model with intermediate vertices in Section~\ref{sec-CM-IV}, and the Preferential Attachment model (PAM) in Section~\ref{sec-PAM}. For the undirected preferential attachment model, we use the basic version with $m=1$ and $\delta=0$, which implies $\gamma=2$. In both configuration models (without and with intermediate vertices) we generate the degree sequences by rounding up i.i.d. values of a continuous random variable $\eta$ with Pareto distribution: $\P(\eta>x)=4x^{-2}$, $x>2$. The exponent $\gamma=\chr{2}$ is chosen for a fair comparison to PAM, and all degrees are at least three for the strongest disassortativity in the model with intermediate in the model with intermediate vertices, see (\ref{eq:spearman_cmie}). In case of the configuration graph in Section~\ref{sec-CM}, we consider two versions: the original model with self-loops and double edges present, and the model where self-loops and double-edges are removed. The rank correlation coefficient $\rho^{\rm rank}(G)$ is computed as in (\ref{eq:rhoG_rank}). The results are presented in Table~\ref{tab:RG}.
\begin{table}[htb]%
{\small\centerline{
\begin{tabular}{|c|l|c|c|c|c|}
\hline
&&\multicolumn{4}{|c|}{$n$}\\
\cline{3-6}
Model&Characteristic&$10^2$&$10^3$&$10^4$&$10^5$\\
\hline
&$\E_N(\rho(G_n))$&-0.0070&-0.0018&-0.0011&0.0006\\
Configuration model&$\sigma_N(\rho(G_n))$&
0.0735&0.0221&0.0077&0.0017\\
\cline{2-6}
with self-loops and double edges&$\E_N(\rho^{\rm rank}(G_n))$&
0.0056&-0.0098&-0.0036&0.0005\\
&$\sigma_N(\rho^{\rm rank}(G_n))$&
0.0504&0.0150&0.0046&0.0019\\
\hline
&$\E_N(\rho(G_n))$&-0.0713&-0.0226&-0.0150&-0.0032\\
Configuration model&$\sigma_N(\rho(G_n))$&0.0546&0.0188&0.0092&0.0029\\
\cline{2-6}
without self-loops and double edges&$\E_N(\rho^{\rm rank}(G_n))$&
-0.0409&-0.0094&-0.0032&-0.0006\\
&$\sigma_N(\rho^{\rm rank}(G_n))$&0.0700&0.0201&0.0083&0.0021\\
\hline
&$\E_N(\rho(\bar{G}_n))$&-0.2804&   -0.1346&   -0.0572&   -0.0291\\
Configuration model&$\sigma_N(\rho(\bar{G}_n))$&0.0742&0.0517&0.0279&0.0147\\
\cline{2-6}
with intermediate vertices&$\displaystyle{\E_N(\rho^{\rm rank}(\bar{G}_n))}$&
-0.7523&-0.7498&-0.7498&-0.7500\\
&$\sigma_N(\rho^{\rm rank}(\bar{G}_n))$&
0.0081&0.0025&0.0008&0.0003\\
\hline
&$\E_N(\rho(G_n))$&-0.2682&-0.1282&-0.0608&-0.0272\\
Preferential attachment&$\sigma_N(\rho(G_n))$&
0.0575&0.0271&0.0132&0.0064\\
\cline{2-6}
&$\E_N(\rho^{\rm rank}(G_n))$&-0.4347&-0.4263&-0.4288&-0.4289\\
&$\sigma_N(\rho^{\rm rank}(G_n))$& 0.0627&0.0272&0.0065&0.0020\\
\hline
\end{tabular}}
\caption{\small Estimated mean and standard deviation of $\rho(G_n)$ and $\rho^{\rm rank}(G_n)$ obtained from 20 realizations of $G_n$ \chn{for random graph models in Sections~\ref{sec-CM}--\ref{sec-PAM}}.}
\label{tab:RG}
}
\end{table}

\vskip0.5cm

The results for the configuration model with intermediate vertices confirm our findings in Section~\ref{sec-CM-IV}: \chn{Pearson's} coefficient converges to zero, while Spearman's rho quickly converges to $-0.75$ revealing the strong negative dependence. For \chrr{the} PAM, Pearson's coefficient 
converges to zero, as indicated in Theorem~\ref{thm-corr-coef-pos-scale-free},
while Spearman's rank correlation clearly indicates a negative dependence. This can be understood by noting that the majority of edges of vertices
with high degrees, which are old vertices,  come from vertices which are added
late in the graph growth process and thus have small degree. On the other hand,
by the growth mechanism of the PAM,  vertices with low degree are
more likely to be connected to vertices having high degree, which indeed suggests
negative degree-degree dependencies. 

We emphasize that under given model assumptions, the graphs of different sizes have been constructed by the same algorithm. Thus, their mixing patterns are exactly the same. As we predicted, the Pearson correlation coefficient fails to reflect the intrinsic properties of the model because its absolute value decreases with the graph size, and converges to zero for all models. On the contrary, Spearman's rho consistently shows neutral mixing for the classical configuration model, moderately disassortative mixing for the Preferential Attachment graph, and strongly disassortative mixing for the configuration model with intermediate vertices.

\subsection{Numerical results for collections of bipartite graphs}
We next compute the  degree-degree dependencies in the collection of bipartite graphs 
discussed in Section \ref{sec:random_assortativity}.
In Table~\ref{tab:bipartite} we present numerical results for 
$\rho(G_n)$ and $\rho^{\rm rank}(G_n)$. Here we choose 
$b=1/2$, $a=2$, $\chrr{\xi}$ has a generalized Pareto distribution 
$\P(\chrr{\xi}>x)=(1+(x-1)/2.8)^{-2.8}$, $x> 1$, and the degrees $X$ and $Y$ are obtained by rounding up the values in \eqref{eq:X}.
	\begin{table}[htb]%
	{\small\centerline{
	\begin{tabular}{|l|c|c|c|c|}
	\hline
	$n$&$10^2$&$10^3$&$10^4$&$10^5$\\
	\hline
	$\E_N(\rho(G_n))$&0.6554&0.7247&0.8042&0.8265\\
	$\sigma_N(\rho(G_n))$&0.1145&0.1406&0.0689&0.0654\\
	\hline
	$\E_N(\rho^{\rm rank}(G_n))$&0.7575&0.7950&0.8526&0.8615\\
	$\sigma_N(\rho^{\rm rank}(G_n))$&
	 0.0735&0.1377&0.0218&0.0074\\
	\hline
	\end{tabular}}}
\caption{\small Estimated mean and standard deviation of $\rho(G_n)$ and $\rho^{\rm rank}(G_n)$ for the collection of $n$ complete bipartite graphs. The number of realizations for 
each graph size is 20.}
\label{tab:bipartite}
\end{table}

\begin{figure}[ht]
	\centerline{\includegraphics[height=3in]{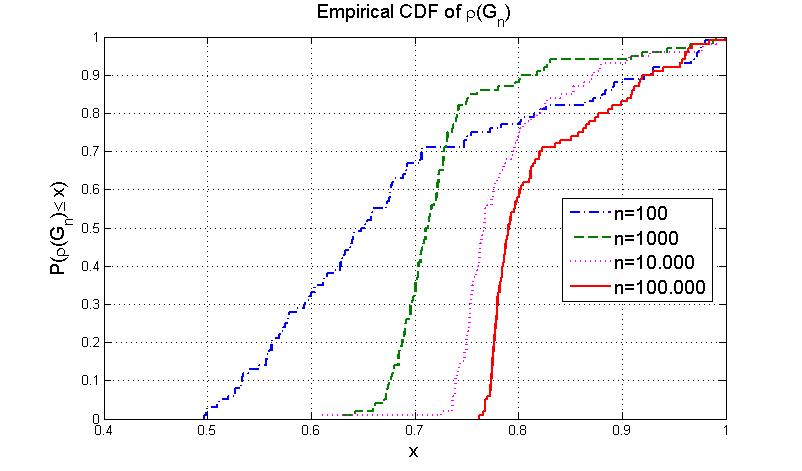}}
	\caption{\small The empirical distribution function $\P(\rho(G_n)\le x)$ for
	$100$ observed values of $\rho(G_n)$, where $G_n$ is a collection of $n$ 
	complete bipartite graphs.}
	\label{fig:bipartite}
	\end{figure}

Note that in this model there is a genuine dependence between the correlation 
measure and the graph size. Indeed, if $n=1$ then the assortativity coefficient 
equals $-1$ because nodes with larger degrees are connected to nodes with 
smaller degrees. However, when the graph size grows, the positive correlations 
start dominating because of the positive linear dependence between $X$ and $Y$. 
We see that again the rank correlation captures the relation faster and gives 
consistent results with decreasing dispersion of values. Finally, Figure~\ref{fig:bipartite} 
shows the changes in the empirical distribution of $\rho(G_n)$ as $n$ grows. 
It is clear that a part of the probability mass is spread over the interval $(0.8,1)$. In the limit, 
$\rho(G_n)$ has a non-zero density on this interval. The difference between the 
crossproducts and the expectation squared in $\rho(G_n)$ is only of the order 
$n^{1-2/\gamma}$, which is about $n^{0.29}$ in our example, thus, the convergence is 
too slow to be observed at $n=100.000$. 

\subsection{Web samples and social networks} For completeness, we present the numerical results for web samples and social networks from \cite{HofLit13}, see in Table~\ref{tab:data}. We used the compressed graph data from the Laboratory of Web Algorithms (LAW) at the Universit\`a degli studi di Milano~\cite{Boldi2004,Boldi2011} with {\rm bvgraph} MATLAB package~\cite{Gleich2010ComputingPageRank}. The {\it stanford-cs} database~\cite{Constantine2007} is a 2001 crawl that includes all pages in the cs.stanford.edu domain. In datasets (iv), (vii), (viii) we evaluate $\rho(G_n)$, $\rho^{\rm rank}(G_n)$ and $\rho^-(G_n)$ (see (\ref{eq:rho-})) over 1000 random edges, and present the average over 10 such evaluations (in 10 samples of 1000 edges, the observed dispersion of the results was small). 

We note that $\rho^{\rm rank}(G_n)$ here is an approximation of (\ref{eq:rhoG_rank}) computed as described in \cite{HofLit13}: we define the random variables $X$ and $Y$ as the degrees on two ends of a random {\it undirected} edge in a graph (that is, here $(u,v)$ and $(v,u)$ represent the same edge); for each edge, when the observed degrees are $a$ and $b$, we assign $[X=a,Y=b]$ or $[X=b,Y=a]$ with probability 1/2; the ties are resolved randomly as in (\ref{eq:rhoG_rank}). The experiments on random graphs show that the values obtained by this algorithm are very close to those computed by (\ref{eq:rhoG_rank}).

\begin{table}{\footnotesize
\begin{tabular}{|c|l|l|c|c|c|c|c|c|}
\hline
nr&Dataset&Description&\# nodes& \# edges& max degree&$\rho(G_n)$&$\rho^{\rm rank}(G_n)$&${\rho}^-(G_n)$\\
\hline
(i)&stanford-cs&web domain&9,914&54,854&340&-.1656&-.1627&-.4648\\\hline
(ii)&eu-2005&.eu web domain&862,664&5,477,938&68,963&-.0562&-.2525&-.0670\\\hline
(iii)&uk@100,000&.uk web crawl&100,000&5,559,150&55,252&-.6536&-.5676&-1.117\\\hline
(iv)&uk@1,000,000&.uk web crawl&1,000,000&77,123,940&403,441&-.0831&-.5620&-.0854\\\hline
(v)&enron&e-mail exchange&69,244&506,898&1,634&-.1599&-.6827&-.1932\\\hline
(vi)&dblp-2010&co-authorship&326,186&1,615,400&238&.3018&.2604&-.7736\\\hline
(vii)&dblp-2011&co-authorship&986,324&6,707,236&979&.0842&.1351&-.2963\\\hline
(viii)&hollywood-2009&co-starring&1,139,905&113,891,327&11,468&.3446&.4689&-0.6737\\\hline
\end{tabular}
}
\caption{(i)--(iv) Web crawls: nodes are web pages, and an (undirected) edge means that there is a hyperlink from one of the two pages to another; (iii),(iv) are  breadth-first crawls around one page. (v) e-mail exchange by Enron employees (mostly part of the senior management): node are employees, and an edge means that an e-mail message was sent from one of the two employees to another. (vi), (vii) scientific collaboration networks extracted from the DBLP bibliography service: each vertex represents a scientist and an edge means a co-authorship of at least one article. (viii) vertices are actors, and two actors are connected by an edge if they appeared in the same movie.}
\label{tab:data}
\end{table}
The most remarkable result here is obtained on the two .uk crawls (iii) and (iv). Here $\rho(G_n)$ is significantly smaller in magnitude on a larger crawl. Intuitively, mixing patterns should not depend on the crawl size. This is indeed confirmed {by} the value of Spearman's {rho}, which consistently shows strong negative correlations in both crawls. We could not observe a similar phenomenon so sharply in (vi) and (vii), probably because a larger co-authorship network incorporates articles from different areas of science, and the culture of scientific collaborations can vary greatly from one research field to another. 

We also notice that, as predicted by our results, the \chrr{small} in magnitude values of ${\rho}^-(G_n)$ result in profound difference in magnitude between $\rho(G_n)$ and $\rho^{\rm rank}(G_n)$.  This is clearly seen in the data sets (ii), (iv) and (v). Again, (ii) and (iv) are the largest among the analyzed web crawls.

The observed behaviour of Pearson's coefficient is explained by the results proved in this paper in that $\rho(G_n)$ is strongly influenced by the large dispersion in the degree values, and particularly by
the presence of hubs. The latter increases with graph size because of the scale-free phenomenon. As a result, $\rho(G_n)$ becomes smaller in magnitude \chrr{when $n$ increases,} which makes it impossible to compare graphs of different sizes. In contrast, the {\it ranks} of the degrees are drawn from a uniform distribution on $[0,1]$, scaled by the factor $|E'|$. Clearly, when a correlation coefficient is computed, the scaling factor cancels, and therefore Spearman's rho provides consistent results in the graphs of different sizes.

\section{Discussion}
\label{sec:discussion}
In this paper, we have investigated dependency measures for power-law random variables. 
We have argued that Pearson's correlation coefficient, despite its appealing feature that it is always in $[-1,1]$, is inappropriate to describe dependencies between heavy-tailed random variables. Indeed, the two main problems with the 
sample correlation coefficient are that (a) it can 
converge to a proper random variable when the sample size tends to infinity, indicating that
it fluctuates tremendously as the sample size increases, and (b) that it is always 
asymptotically non-negative when
dealing with non-negative random variables (even when these are obviously negatively 
dependent). In the context of random graphs, the first deficiency means that 
Pearson's coefficient can have a non-vanishing variance even when the size of the graph
is huge, the second \chn{mistakenly suggests} that there do not exist asymptotically disassortative scale-free
graphs. We give proofs for the facts stated above, and illustrate the results using simulations.

Rank correlations are a special case of the broader concept of copulas that are widely used in multivariate analysis, in particular in applications in mathematical finance and risk management. There is a heated discussion in this area about the adequacy and informativeness of such measures, see e.g.\ \cite{Mikosch2006copulas} and consequent reactions. There are several points of criticism. In particular, Spearman's rho uses rank transformation, which changes the observed values of the degrees. Then, first of all, what exactly does Spearman's rho tell us about the dependence between the original values? Second of all, no substantial justification exists for the rank transformation, besides its mathematical convenience. 
We thus do not claim that Spearman's rho is {\it the} solution to the problem. {Nevertheless, compared to the Pearson's coefficient, Spearman's rho has a significant advantage that it is free from the undesirable size-dependency, and converges to meaningful value in the infinite volume limit.}

\chn{We note \chrr{that} Spearman's rho has computational complexity $O(n\log(n))$ because the values of the random variables must be ranked first. Pearson's correlation coefficient is easier to evaluate because it uses the values of the degrees directly, and has computational complexity $O(n)$. Efficient methods for computing Spearman's rho in large graphs is an interesting topic for future research.}

Raising the discussion to a higher level, random variables $X$ and $Y$ are positively 
dependent when a large realization of $X$ typically implies a large realization of $Y$.
A strong form of this notion is when $\P(X>x, Y>y)\geq \P(X>x)\P(Y>y)$
for every $x,y\in {\mathbb R}$, but for many purposes this notion is too
restrictive.  The covariance for non-negative random variables
is obtained by integrating the above inequality over $x,y\geq 0$, so that
it is true for `typical' values of $x,y$. In many cases, however, we
are particularly interested in certain values of $x,y$. 
Another class of methods for measuring rank correlations is
based on the  \emph{angular measure}, a notion originating 
in the theory of multivariate extremes, for which the above inequality
is investigated for \emph{large} $x$ and $y$, so that it describes 
the {\it tail dependence} for a random vector $(X,Y)$, that is, the dependence
between extremely large values of $X$ and $Y$, see e.g.\ \cite{Resnick_HT}. Such tail
dependence is characterized by an probability-like measure, or, the {angular measure}, on $[0,1]$. Informally, a
concentration of the angular measure around the points 0 and 1 indicates
independence of large values, while concentration around some other number
$a\in(0,1)$ suggests that a certain fraction of large values of
$Y$ comes together with large values of $X$.
In~\cite{Volkovich08WWW,Volkovich2009complex} a first attempt
was made to compute the angular measure between 
in-degree of a node and its importance measured by the
Google PageRank algorithm. Strikingly, completely different dependence structures were
discovered in Wikipedia (independence), Preferential Attachment
networks (complete dependence) and the Web (intermediate case).



\paragraph{Acknowledgments.}
We thank Yana Volkovich for the code generating a preferential attachment graph \chr{and Marie Albenque 
for a counter example that shows that negative dependence of $(X,Y)$ in does not follow from negative dependence of $(X_1,Y_1)$.} \chr{We further thank Juli Komj\'athy for her careful reading of the paper, which has tremendously improved the presentation and has corrected several typos and errors.}
{This article is also the result of joint research in the 3TU
Centre of Competence NIRICT (Netherlands Institute for Research on
ICT) within the Federation of Three Universities of Technology in
The Netherlands. The work of RvdH was supported in part by the Netherlands Organisation for Scientific Research (NWO). The work of NL is partially supported by the EU-FET Open grant NADINE (288956).}

\bibliographystyle{plain}
\bibliography{myrefs}

\end{document}